\documentclass[11pt]{article}
\usepackage{epsfig}
\usepackage{amsthm}
\usepackage{amssymb}
\usepackage{amsmath}
\usepackage{amsfonts}
\usepackage{color}
\newcommand{\field}[1]{\mathbb{#1}}
\newcommand{\R}{\field{R}}

\newcommand{\N}{\field{N}}
\newcommand{\Z}{\field{Z}}

\newcommand{\E}{\field{E}}
\def\qed{\hfill$\diamondsuit$}

\theoremstyle{example} \theoremstyle{remark} \theoremstyle{lemma}
\theoremstyle{definition} \theoremstyle{corol}
\theoremstyle{proposition} \theoremstyle{condition}
\theoremstyle{assumption}
\newtheorem{assumption}{\n{Assumption}}[section]
\newtheorem{theorem}{\n{Theorem}}[section]

\newtheorem{remark}{\n{Remark}}[section]
\newtheorem{lemma}{\n{Lemma}}[section]

\font\n=cmcsc10
\def\cov{{\mbox{cov}}}
\def\var{{\mbox{var}}}
\def\cum{{\mbox{cum}}}
\begin{document}

\bibliographystyle{plain}

\centerline{\Large A GENERALIZED PORTMANTEAU TEST }

\centerline{\Large FOR INDEPENDENCE BETWEEN }

\centerline{\Large TWO STATIONARY TIME SERIES \footnote{ I am grateful to the coeditor Robert Taylor and
the referees for useful suggestions  that greatly improved the paper. I also want
 to thank Hernando Ombao and Wei Biao Wu for helpful comments on an earlier version. The research is supported in 
 part by NSF grant DMS-0804937. Address correspondence to: Xiaofeng Shao, Department of
Statistics, University of Illinois at Urbana-Champaign,  725 South
Wright St, Champaign, IL, 61820; e-mail: xshao@uiuc.edu}}

\bigskip
 \centerline{\textsc{By Xiaofeng Shao}}
\centerline{\today} \centerline {\it University of Illinois at
Urbana-Champaign}

\bigskip

\noindent We propose generalized portmanteau-type test statistics
in the frequency domain to test independence between two
stationary time series. The test statistics are formed analogous
to the one in Chen and Deo (2004, {\it Econometric Theory} 20,
382-416), who extended the applicability of portmanteau
goodness-of-fit test to the long memory case.
Under the null hypothesis of independence, the asymptotic standard
normal distributions of the proposed statistics  are derived under
fairly mild conditions. In particular, each time series is allowed
to possess short memory,  long memory or anti-persistence. A
simulation study shows that the tests have reasonable size and power
properties.

 \pagenumbering{arabic}

\setcounter{page}{1}
\section{Introduction}

Recently there has been considerable amount of work devoted to
testing the independence or non-correlation of two stationary time
series, e.g. Bouhaddioui and Roy (2006), Duchesne and Roy (2003),
Eichler (2006), Hallin and Saidi (2005, 2007). Consider two univariate
time series
\begin{eqnarray}
\label{eq:xtyt}
 X_{1t}=\sum_{j=0}^{\infty}a_ju_{t-j},~~~X_{2t}=\sum_{j=0}^{\infty} b_j
 v_{t-j},~t\in\Z,
 \end{eqnarray}
where $\{u_t\}$ and $\{v_t\}$ are each a sequence of independent
and identically distributed (iid) random variables. Our goal is to
test the null hypothesis that $\{X_{1t}\}$ and $\{X_{2t}\}$ are
uncorrelated at all lags. If the joint distribution of the two
series is Gaussian, then non-correlation is equivalent to
independence. Given a realization $\{X_{kt}\}_{t=1}^{n}, k=1,2$,
of length $n$ of the processes defined in (\ref{eq:xtyt}),
a popular portmanteau-type statistic is based on the sum of the
weighted squared cross-correlations, i.e.
\begin{eqnarray}
\label{eq:Gn} G_n=\sum_{j=1-n}^{n-1} K^2(j/B_n)
\hat{\rho}_{\hat{u}\hat{v}}^2(j),
\end{eqnarray}
 where $K(\cdot)$
is a kernel function, $B_n$ is the bandwidth and
$\hat{\rho}_{\hat{u}\hat{v}}(j)$ is the (empirical)
cross-correlation of two residual series $\{\hat{u}_t\}$ and
$\{\hat{v}_t\}$ at lag $j$.  Typically the residual series are
obtained by prewhitening the two time series separately, which can
be done by fitting a parametric model, such as an ARMA model to
each time series [Haugh (1976)]. A multivariate extension of
Haugh's (1976) idea  can be found in El Himdi and Roy (1997) and
Pham et al. (2003). However, the asymptotic null distribution of
the test statistic is invalid if the parametric models are
misspecified. To avoid possible model misspecifications, a long
autoregression can be fitted to each time series and the
consistency of the test is ensured as long as the orders of
autoregressions grow properly as the sample size increases [Hong
(1996)]. See Bouhaddioui and Roy (2006) for a multivariate
extension of Hong's (1996) test statistic.

All of the work mentioned above seems to exclude long memory time
series.
  In the past two decades, a great deal of research  has been conducted on long memory
time series and its importance has been recognized in both applied
and theoretical time series literature [Doukhan et al. (2003),
Robinson (2003), Teyssi$\grave{e}$re and Kirman (2005)]. In this
article, we develop tests for independence between two long memory
time series, which seem lacking in the literature. Our test
statistics are frequency domain analogues of Hong's (1996) test
statistic, which
 has been  recently reformulated by Eichler (2006) in the frequency domain.
However, the theory in the latter paper was developed under
stringent conditions. For example, the author assumed the
existence of all moments and summability of joint cumulants up to
all orders, which ruled out the interesting long memory case.

The portmanteau-type statistic for testing the independence
between two time series bears some resemblance to that for testing
the goodness-of-fit for time series models. For the latter
problem, Chen and Deo (2004a) formulated a generalized test
statistic in the frequency domain and extended the applicability
of portmanteau-type test to the long memory case. The development
in this paper parallels  Chen and Deo (2004a) in that for each time
series the prewhitening is also done in the frequency domain with
the unknown spectral density replaced by the estimated spectral
density. However, our work differs from Chen and Deo (2004a) in
two important aspects. First, in our feasible test statistics,
the spectral densities of $X_{kt}$, $k=1,2$, can be estimated
under parametric assumptions [see Chen and Deo (2004a)], or via a
nonparametric approach, where
 no parametric models need to be specified. See
Sections~\ref{subsec:1} and \ref{subsec:2} for more details.
 Second, our results are applicable to not only short/long memory
time series, but also the anti-persistent case. Similar to correlation-based independence tests, the tests proposed can't
discriminate between non-correlation and independence for
non-Gaussian time series, although our asymptotic theory allows for non-Gaussian linear processes.
Nevertheless, they still help the practitioner extract
useful information from the data, especially when the null
hypothesis is rejected. Finally, we note that the techniques and results developed here are not directly applicable  to relax
 the Gaussian assumption in Chen and Deo (2004a), where the proofs require some sharp bounds for the products of discrete Fourier transforms.





Now we introduce some notation. For a vector ${\pmb{x}} = (x_1,
\cdots, x_q)'\in \R^q$, let $|\pmb{x}| = (\sum_{i=1}^q
x_i^2)^{1/2}$. For two sequences $(a_n)$, $(b_n)$, denote
by $a_n \sim b_n$ if $a_n / b_n \to 1$ as $n \to \infty$. Denote
by $\rightarrow_{D}$ and ${\rightarrow}_{p}$ convergence in
distribution and in probability, respectively. The symbols
$O_{p}(1)$ and $o_{p}(1)$ signify being bounded in probability and
convergence to zero in probability respectively. Let $N(0,1)$ be
the standard normal distribution.

The paper is organized as follows. In Section~\ref{sec:2}, we
introduce our test statistics and derive their asymptotic null
distributions. The size and power properties of our tests are
examined in Section~\ref{sec:3} through simulations. The technical
details are relegated to the Appendix.

\section{The Test Statistics and Their Asymptotic Null Distributions}
\label{sec:2}
 Throughout the paper, we assume $X_{kt}$, $k=1,2$,
admit the linear processes of the form (\ref{eq:xtyt}),
where the innovations $u_t$ and $v_t$ satisfy the following
assumptions.
\begin{assumption}
\label{as:innovation}{\rm The series $\{u_t\}$ and $\{v_t\}$ are
each a sequence of iid random variables with mean zero and finite
fourth moment. Without loss of generality, we assume
$\var(u_t)=\var(v_t)=1$. Denote by $c_4(u)=\cum(u_0,u_0,u_0,u_0)$
and $c_4(v)=\cum(v_0,v_0,v_0,v_0)$.}
\end{assumption}

Let $i=\sqrt{-1}$ be the imaginary unit. For a complex number $c$
let $\overline{c}$ be its conjugate. For any two processes
$\{Z_{1t},Z_{2t},t\in\Z\}$, denote by $f_{Z_1Z_1}(\lambda)$ and
$f_{Z_2Z_2}(\lambda)$ their spectral densities respectively;
define their Fourier transforms, periodograms and cross
periodogram by
\[w_{Z_k}(\lambda)=\frac{1}{\sqrt{2\pi n}}\sum_{t=1}^{n}Z_{kt} e^{it\lambda}, ~I_{Z_k}(\lambda)=|w_{Z_k}(\lambda)|^2,~\mbox{and}~I_{Z_1 Z_2}(\lambda)=w_{Z_1}(\lambda)\overline{w_{Z_2}(\lambda)}\]
for $k=1,2$.  Let $\lambda_j=2\pi j/n, j=1,2,\cdots,n$, be the
Fourier frequencies. To motivate our test statistic,  we express
(\ref{eq:Gn}) in the frequency domain. Define the residual
cross-covariance function by
\[\hat{R}_{\hat{u}\hat{v}}(j)=\left\{\begin{array}{cc}
n^{-1}\sum_{t=j+1}^{n} \hat{u}_{t-j}\hat{v}_{t}&(j\ge 0)\\
n^{-1}\sum_{t=1-j}^{n}\hat{u}_{t}\hat{v}_{j+t}&(j<0)\end{array}\right.\]
and the auto-covariance function for the residual series
$\hat{u}_t$ and $\hat{v}_t$ by
\[\hat{R}_{\hat{u}\hat{u}}(j)=n^{-1}\sum_{t=1+|j|}^{n}\hat{u}_t\hat{u}_{t-|j|},~~\hat{R}_{\hat{v}\hat{v}}(j)=n^{-1}\sum_{t=1+|j|}^{n}\hat{v}_t\hat{v}_{t-|j|}.\]
 Then
$\hat{\rho}_{\hat{u}\hat{v}}(j)=\hat{R}_{\hat{u}\hat{v}}(j)/\{\hat{R}_{\hat{u}\hat{u}}(0)\hat{R}_{\hat{v}\hat{v}}(0)\}^{1/2}$.
 Let
\[\hat{f}_{Z_1 Z_2}(\lambda)=(2\pi)^{-1}\sum_{|j|<n}K(j/B_n)\hat{R}_{Z_1 Z_2}(j)
e^{-i\lambda j},~Z_1,Z_2=\hat{u},\hat{v},\]
 where $K(\cdot)$ is a symmetric kernel function
with $K(0)=1$. Then
\begin{eqnarray}
\label{eq:Gn2}
G_n=2\pi\int_{0}^{2\pi}|\hat{f}_{\hat{u}\hat{v}}(\lambda)|^2
d\lambda
\left(\int_{0}^{2\pi}\hat{f}_{\hat{u}\hat{u}}(\lambda)d\lambda
\int_{0}^{2\pi}\hat{f}_{\hat{v}\hat{v}}(\lambda)d\lambda\right)^{-1}.
\end{eqnarray}
 Let $W(\lambda)=(2\pi)^{-1}\sum_{|h|<n} K(h/B_n) e^{-i h\lambda}$
be the spectral window function corresponding to the kernel
function $K(\cdot)$. We
  can write $\hat{f}_{\hat{u}\hat{v}}(\lambda)$,  $\hat{f}_{\hat{u}\hat{u}}(\lambda)$ and $\hat{f}_{\hat{v}\hat{v}}(\lambda)$ into  the following
  equivalent forms in the frequency domain (see Equation (3) in Chen and Deo (2004a)),
  \begin{eqnarray}
  \label{eq:fuv}
  \hat{f}_{Z_1Z_2}(\lambda)=\int_{0}^{2\pi}W(\lambda-w)I_{Z_1Z_2}(w)dw,
  ~Z_1,Z_2=\hat{u}, \hat{v}.
  \end{eqnarray}

 The  expressions (\ref{eq:Gn2}) and  (\ref{eq:fuv})
  motivate us to propose the following test statistic:
  \begin{eqnarray*}
  \label{eq:Tn}
T_n&=&\frac{4\pi^2}{n}\sum_{l=0}^{n-1}|\hat{f}_{X_1X_2}(\lambda_l)|^2
\left\{\frac{2\pi}{n}\sum_{l=0}^{n-1}\hat{f}_{X_1X_1}(\lambda_l)
\frac{2\pi}{n}\sum_{l=0}^{n-1}\hat{f}_{X_2X_2}(\lambda_l)\right\}^{-1},
\end{eqnarray*}
where
\begin{eqnarray*}
\hat{f}_{X_h
X_k}(\lambda)=\frac{2\pi}{n}\sum_{j=1}^{n-1}\frac{W(\lambda-\lambda_j)I_{X_h
X_k}(\lambda_j)}{\sqrt{f_{X_hX_h}(\lambda_j)f_{X_kX_k}(\lambda_j)}},
~h,k=1,2.
\end{eqnarray*}
Similar to Chen and Deo (2004a), our test statistic $T_n$ is
obtained by discretizing the integrals in (\ref{eq:Gn2}) and
(\ref{eq:fuv}) with $I_{\hat{u}\hat{v}}$ replaced by
$I_{X_1X_2}/\sqrt{f_{X_1X_1}f_{X_2X_2}}$ and $I_{\hat{u}\hat{u}}$
($I_{\hat{v}\hat{v}}$) replaced by $I_{X_1X_1}/f_{X_1X_1}$
($I_{X_2X_2}/f_{X_2X_2}$). Since $\hat{f}_{X_1X_2}$,
$\hat{f}_{X_1X_1}$ and $\hat{f}_{X_2X_2}$ are evaluated only at
Fourier frequencies, our test statistic $T_n$ is mean-invariant.
Our test statistic $T_n$ is infeasible since we do not know the
true spectral densities $f_{X_1X_1}(\lambda)$ and
$f_{X_2X_2}(\lambda)$ in practice. In the next two subsections, we
introduce two ways to estimate the spectral densities, which lead
to two feasible test statistics.


To establish the asymptotic null distribution of $T_n$, we make
the following assumptions on the kernel function $K(\cdot)$ and
the bandwidth $B_n$.

\begin{assumption}
\label{as:kernel} {\rm The kernel function $K:\R\rightarrow
[-1,1]$ has compact support on $[-1,1]$. It is  differentiable
except at a finite number of points and symmetric with $K(0)=1$.}
\end{assumption}
The assumption that the kernel function has compact support can be
relaxed; see Chen and Deo (2004a). Here we decide to retain this
assumption to avoid
 more technical complications in view of our long and technical proof.  It is worth noting that several
commonly-used kernels in spectral analysis, such as Bartlett,
Parzen and Tukey kernels,  satisfy Assumption~\ref{as:kernel} (see
Priestley (1981), p 446-447).

\begin{assumption}
\label{as:bandwidth} {\rm The bandwidth $B_n$ satisfies $(\log^2
n)/B_n\rightarrow 0$ and $B_n(\log^2 n)/n\rightarrow 0$.}
\end{assumption}

Let $A(\lambda)=\sum_{j=0}^{\infty}a_j e^{ij\lambda}$ and
$B(\lambda)=\sum_{j=0}^{\infty}b_je^{ij\lambda}$. The following
assumption is made regarding the long memory behavior of
$\{X_{kt}\}$, $k=1,2$ and is satisfied by two commonly-used long
memory time series models: FARIMA (fractional autoregressive
integrated moving average) models and fractional Gaussian noise
[Beran (1994)].

\begin{assumption}
\label{as:dependence} {\rm For $k=1,2$, assume
$f_{X_kX_k}(\lambda)\sim |\lambda|^{-2d_{k0}}G_k$ as
$\lambda\rightarrow 0$, where $d_{k0}\in (-1/2,1/2)$ and $G_k\in
(0,\infty)$. Further we assume that
\[|\partial A(\lambda)/\partial \lambda|=O(|A(\lambda)||\lambda|^{-1})~~\mbox{and}~~|\partial B(\lambda)/\partial \lambda|=O(|B(\lambda)||\lambda|^{-1}) \] hold uniformly in $\lambda\in (0,\pi]$.
}
\end{assumption}
For $k=1,2$, the process $\{X_{kt}\}$ is said to possess long
memory if $d_{k0}\in (0,1/2)$, short memory if $d_{k0}=0$ and
anti-persistence if $d_{k0}\in (-1/2,0)$. Our results cover the
short memory and anti-persistent cases as well.

\begin{theorem}
\label{th:mainresult1} Suppose that the two processes $\{u_t\}$
and $\{v_t\}$ are independent. Under
Assumptions~\ref{as:innovation}-\ref{as:dependence}, we have
\[\frac{n T_n- B_n s(K)}{\sqrt{2 B_n d(K)}}\rightarrow_{D} N(0,1),\]
where $s(K)=\int_{-\infty}^{\infty}K^2(x)dx$ and
$d(K)=\int_{-\infty}^{\infty}K^4(x)dx$.
\end{theorem}


\subsection{Feasible Test Statistic I}
\label{subsec:1}

In this subsection, we estimate the spectral densities in a
parametric way, so the following
 parametric assumptions are imposed on $f_{X_kX_k}(\lambda)$, $k=1,2$.


\begin{assumption}
\label{as:parametric} {\rm For $k=1,2$, let the spectral density
of the process $\{X_{kt}\}$ be
$f_k(\lambda;{\bf{{\pmb{\theta}}}}_{k0})$, where
${\pmb{\theta}}_{k0}$ is the true parameter vector that lies in
the interior of the compact set $\Theta_{k}\subset \R^{q_k}$,
$q_k\in\N$.
Suppose that the estimator $\hat{{\pmb{\theta}}}_k$ satisfies
$|\hat{{\pmb{\theta}}}_k-{\pmb{\theta}}_{k0}|=O_{p}(n^{-1/2})$,
$k=1,2$.}
\end{assumption}
Here  we can take  the  Whittle pseudo-maximum likelihood
estimator as $\hat{{\pmb{\theta}}}_k$. The root-$n$ asymptotic
normality of Whittle estimator has been established by Hannan
(1973) for $d=0$, and by Fox and Taqqu (1986), Dahlhaus (1989) and
Giraitis and Surgalis (1990) for $d\in (0,1/2)$. See Velasco and
Robinson (2000) for the case $d\in (-1/2,0)$.

 Let
${\pmb{\theta}}=({\pmb{\theta}}_1',{\pmb{\theta}}_2')'$ and
$\hat{{\pmb{\theta}}}=(\hat{\pmb{\theta}}_1',\hat{{\pmb{\theta}}}_2')'$.
Then we can replace $f_{X_1X_1}(\lambda)$ and
$f_{X_2X_2}(\lambda)$ in $T_n$ with
$f_1(\lambda;\hat{{\pmb{\theta}}}_1)$ and
$f_2(\lambda;\hat{{\pmb{\theta}}}_2)$ respectively and get the
following feasible estimator
\begin{eqnarray*}
T_n(\hat{{\pmb{\theta}}})=\frac{4\pi^2}{n}\sum_{l=0}^{n-1}|\tilde{f}_{X_1X_2}(\lambda_l)|^2\left\{\frac{2\pi}{n}\sum_{l=0}^{n-1}
\tilde{f}_{X_1X_1}(\lambda_l)\frac{2\pi}{n}\sum_{l=0}^{n-1}
\tilde{f}_{X_2X_2}(\lambda_l)\right\}^{-1},
\end{eqnarray*}
\begin{eqnarray*}
~\mbox{where}~~~~\tilde{f}_{X_hX_k}(\lambda)=\frac{2\pi}{n}\sum_{j=1}^{n-1}\frac{W(\lambda-\lambda_j)I_{X_hX_k}(\lambda_j)}{\sqrt{f_h(\lambda_j;\hat{{\pmb{\theta}}}_h)f_k(\lambda_j;\hat{{\pmb{\theta}}}_k)}},~h,k=1,2.
\end{eqnarray*}

To obtain the asymptotic null distribution of
$T_n(\hat{\pmb{\theta}})$, we further make the following two
assumptions; compare Assumptions 6\&7 in Chen and Deo (2004a).

\begin{assumption}
\label{as:longmemo}
{\rm
 For $k=1,2$, let ${\pmb{\theta}}_{k0}=({\pmb{\beta}}_{k0}',d_{k0})'$. The spectral density
$f_k(\lambda;{\pmb{\theta}}_{k0})=f_k^*(\lambda;d_{k0})g_k^*(\lambda;{\pmb{\beta}}_{k0})$,
where $f_k^*$ and $g_k^*$ are even functions on $[-\pi,\pi]$,
$f_k^*(\lambda;d_{k0})\sim c(d_{k0})\lambda^{-2d_{k0}}$
 as $\lambda\rightarrow 0$ for some constant $c(d_{k0})>0$, $g_k^*(\lambda;{\pmb{\beta}}_{k0})$
is bounded away from zero and differentiable on $[-\pi,\pi]$.
Further assume that the $q_k$th component of $\Theta_{k}$ is
contained in the closed interval $[-\kappa,\kappa]$ for some
$\kappa\in (0,0.5)$. }
\end{assumption}

\begin{assumption}
\label{as:condiff} {\rm For $k=1,2$,
\begin{enumerate}
 \item  $\log f_k(\lambda;{\pmb{\theta}}_k)$,  $\partial\log
f_k(\lambda;{\pmb{\theta}}_k)/\partial {\pmb{\theta}}_{k_u}$ and
$\partial^2 \log
f_k(\lambda;{\pmb{\theta}}_k)/\partial\theta_{k_u}\partial
\theta_{k_v}$ are continuous at all $(\lambda,{\pmb{\theta}}_k)$
except $\lambda=0$. Further,
\[\sup_{\lambda\in [0,2\pi]}\sup_{{\pmb{\theta}}_k\in\Theta_{k}}|\lambda|^{2d_{k0}}f_k(\lambda;{\pmb{\theta}}_{k})=A_k,~\mbox{for
some}~A_k\in(0,\infty).
\]
\item For any $\delta>0$,
\begin{eqnarray}
\label{eq:111stder} \sup_{\lambda\in
[0,2\pi]}\sup_{{\pmb{\theta}}_k\in\Theta_{k}}
|\lambda|^{\delta}\left|\frac{\partial \log
f_k(\lambda;{\pmb{\theta}}_k)}{\partial\theta_{k_u}}\right|=A_k,~\mbox{for
some}~A_k\in (0,\infty) \end{eqnarray} and
\begin{eqnarray}
\label{eq:222ndder} \sup_{\lambda\in
[0,2\pi]}\sup_{{\pmb{\theta}}_k\in\Theta_{k}}
|\lambda|^{\delta}\left|\frac{\partial^2 \log
f_k(\lambda;{\pmb{\theta}}_k)}{\partial\theta_{k_u}\partial\theta_{k_v}}\right|=A_k,~\mbox{for
some}~A_k\in (0,\infty).
\end{eqnarray}

\item There exists a constant $C$ such that
\begin{eqnarray}
\label{eq:lipschitz}
|{f_k^{1/2}(\lambda;{\pmb{\theta}}_k^{(1)})}-{f_k^{1/2}(\lambda;{\pmb{\theta}}_k^{(2)})}|\le
C|{\pmb{\theta}}_k^{(1)}-{\pmb{\theta}}_k^{(2)}|{f_k^{1/2}(\lambda;{\pmb{\theta}}_k^{(2)})}
\end{eqnarray}
 uniformly in $\lambda$ and all
${\pmb{\theta}}_k^{(1)}=({\pmb{\beta}}_k^{(1)'},d_k^{(1)})'$ and
${\pmb{\theta}}_k^{(2)}=({\pmb{\beta}}_k^{(2)'},d_k^{(2)})'$,
$d_k^{(1)}<d_k^{(2)}$.

\end{enumerate}
}
\end{assumption}
In Chen and Deo (2004a), the same (and equivalent) form of
conditions as (\ref{eq:111stder}) and (\ref{eq:222ndder}) were
imposed on $\partial f_k^{-1}(\lambda;\pmb{\theta}_{k})/\partial
\theta_{k_u}$ and $\partial^2
f_k^{-1}(\lambda;\pmb{\theta}_k)/\partial \theta_{k_u}\partial
\theta_{k_v}$. The condition~(\ref{eq:lipschitz}) is a variant of
the  Assumption 7(iii) in their paper. It is not hard to verify
Assumptions~\ref{as:longmemo} and \ref{as:condiff} for fractional
Gaussian noise and FARIMA processes.

\begin{theorem}
\label{th:mainresult2}
Suppose that the two processes $\{u_t\}$ and
$\{v_t\}$ are independent.
 Under Assumptions~\ref{as:innovation}-\ref{as:condiff}, we have
\[\frac{n T_n(\hat{{\pmb{\theta}}})- B_n s(K)}{\sqrt{2 B_n d(K)}}\rightarrow_{D} N(0,1).\]
\end{theorem}

 In forming our feasible test statistic
$T_n(\hat{\pmb{\theta}})$, we assume correct parametric
specifications for the spectral densities of $\{X_{kt}\}, k=1,2$.
In practice, model selection technique [Beran et al. (1998)] can
be used to identify the best parametric model for each time
series.

\subsection{Feasible Test Statistic II}
\label{subsec:2}

In this subsection, we shall restrict our attention to the
fractionally differenced autoregressive models of order $\infty$,
denoted as $\mbox{FAR}(\infty,d)$ (i.e. the $d$-th fractional
difference follows an infinite order autoregressive process). Our
consideration of this class of models is motivated by the results
in Bhansali et al. (2006), who obtained the asymptotic normality
of estimated coefficients when fitting an $\mbox{FAR}(p,d)$ model
(i.e. the $d$-th fractional difference follows an autoregressive
model of order $p$) to the observations from an
$\mbox{FAR}(\infty,d)$ process with $p\rightarrow\infty$ as
$n\rightarrow\infty$.

For $k=1,2$, let
\begin{eqnarray}
\label{eq:far}
 (1-B)^{d_{k0}}X_{kt}=Y_{kt},
~~\sum_{j=0}^{\infty}a_{kj}^0Y_{kt}=\varepsilon_{kt},
\end{eqnarray}
where $B$ is the backward shift operator, $\varepsilon_{kt}$ are
iid mean zero random variables with
$\sigma_{k}^2=\E\varepsilon_{kt}^2$ and
$\E\varepsilon_{kt}^4<\infty$.  Assume $a_{k0}^{0}=1$ and
$\sum_{j=0}^{\infty}|a_{kj}^0|<\infty$. Then $X_{kt}$ has the
spectral density
\[f_{X_kX_k}(\lambda)=(2\pi)^{-1}\sigma_{k}^2|1-e^{i\lambda}|^{-2d_{k0}}\left|\sum_{j=0}^{\infty}a_{kj}^0e^{ij\lambda}\right|^{-2}, ~k=1,2.\]
It is easy to see that the processes $\{X_{kt}\}$, $k=1,2$ defined
in (\ref{eq:far}) are a subclass of (\ref{eq:xtyt}). The following
assumption was also made in Bhansali et al. (2006).
\begin{assumption}
\label{as:expondecay} For some $\epsilon>0$,
\[a_k(x):=\sum_{j=0}^{\infty}a_{kj}^0 x^j\not=0, ~|x|<1+\epsilon, ~k=1,2.\]
\end{assumption}

\begin{remark}
\label{rem:akj} {\rm As stated in Remark 2.1 of Bhansali et al.
(2006),  Assumption~\ref{as:expondecay} implies that
$c_k(x)=1/a_k(x)$ has the expansion
\[c_k(x)=\sum_{j=0}^{\infty}c_{kj}^0 x^j,~~~|x|\le 1, ~c_{k0}^{0}=1,~k=1,2\]
and that $\max(a_{kj}^0, c_{kj}^0)\le C r_k^{j},~j\in\N$ for some
constant $C$ and $r_k\in (0,1)$. It is satisfied by FARIMA
processes with finite autoregressive and moving average orders.}
\end{remark}
Following Bhansali et al. (2006), for $k=1,2$, we fit an
$\mbox{FAR}(p_k,d)$ model to the observations
$\{X_{kt},t=1,\cdots,n\}$ generated from the process
(\ref{eq:far}).
Let
$\hat{{\pmb{\gamma}}}_k=(\hat{d}_k,\hat{a}_{k1},\hat{a}_{k2},\cdots,\hat{a}_{kp_k})'$
 be the resulting
estimates of the parameter vector
${\pmb{\gamma}}_k=(d_{k},a_{k1},a_{k2},\cdots,a_{kp_k})'$ in the
$\mbox{FAR}(p_k,d)$ model.
 Note that the spectral density of an $\mbox{FAR}(p_k,d)$ model is
\[f_{k}(\lambda;{\pmb{\gamma}}_k)=\frac{\sigma_k^2}{2\pi}|1-e^{i\lambda}|^{-2d_{k}}\left|\sum_{j=0}^{p_k} a_{kj} e^{ij\lambda}\right|^{-2},~\mbox{where}~a_{k0}=1.\]
 Let ${\pmb{\gamma}}=({\pmb{\gamma}}_1',{\pmb{\gamma}}_2')'$ and $\hat{{\pmb{\gamma}}}=(\hat{{\pmb{\gamma}}}_1',\hat{{\pmb{\gamma}}}_2')'$.
 We form our
test statistic $T_n(\hat{{\pmb{\gamma}}})$ by replacing
$f_{X_kX_k}(\lambda)$ in $T_n$ with
${f}_{k}(\lambda;\hat{{\pmb{\gamma}}}_k)$, i.e.
\[T_n(\hat{{\pmb{\gamma}}})=\frac{4\pi^2}{n}\sum_{l=0}^{n-1}|\check{f}_{X_1X_2}(\lambda_l)|^2\left\{\frac{2\pi}{n}\sum_{l=0}^{n-1} \check{f}_{X_1X_1}(\lambda_l)\frac{2\pi}{n}\sum_{l=0}^{n-1} \check{f}_{X_2X_2}(\lambda_l) \right\}^{-1},\]
where
\[\check{f}_{X_hX_k}(\lambda)=\frac{2\pi}{n}\sum_{j=1}^{n-1}\frac{W(\lambda-\lambda_j)I_{X_hX_k}(\lambda_j)}{\sqrt{{f}_{h}(\lambda_j;\hat{{\pmb{\gamma}}}_h){f}_{k}(\lambda_j;\hat{{\pmb{\gamma}}}_k)}},~h,k=1,2.\]
Since $T_n(\hat{{\pmb{\gamma}}})$ is free of $\sigma_k^2$, we set
$\sigma_k^2=1$, $k=1,2$.

Let ${\pmb{\gamma}}_{k0}=(d_{k0},a_{k1}^0,\cdots,a_{k p_k}^0)'$,
$k=1,2$. We require certain convergence rate for $\hat{{\pmb{\gamma}}}_k$, which has been
 obtained in Bhansali et al's Theorem 3.1.
\begin{assumption}
\label{as:consistent}
 Assume that
$|\hat{{\pmb{\gamma}}}_k-{\pmb{\gamma}}_{k0}|=O_p(\sqrt{p_k/n})$,
$k=1,2$.
\end{assumption}
\begin{theorem}
\label{th:mainresult3}
 Let Assumptions~\ref{as:innovation}-\ref{as:dependence} and \ref{as:expondecay}-\ref{as:consistent} hold. Suppose
 that the two processes $\{\varepsilon_{1t}\}$ and
 $\{\varepsilon_{2t}\}$ are independent. Assume that for $k=1,2$,
 \begin{eqnarray}
 \label{eq:pB1}
 \sum_{l=p_k+1}^{\infty}|a_{kl}^0|=o(\sqrt{B_n}/(n\log
 n)), ~~~p_k\rightarrow\infty,
 \end{eqnarray}
 \begin{eqnarray}
 \label{eq:pB2}
~\mbox{and}~~~p_k^2\log^2 n=o(B_n), ~p_k^2 B_n=o(n).
 \end{eqnarray}
 Then we have
\[\frac{nT_n(\hat{{\pmb{\gamma}}})-B_ns(K)}{\sqrt{2B_n d(K)}}\rightarrow_{D} N(0,1).\]
\end{theorem}
 Let $p_k=l_k\log n$, where
$l_k>-(\log r_k)^{-1}$. Then (\ref{eq:pB1}) holds under
Assumption~\ref{as:expondecay} (see Remark~\ref{rem:akj}) and
(\ref{eq:pB2}) reduces to $\log^4 n=o(B_n)$ and $(\log^2
n)B_n=o(n)$.

\begin{remark}
{\rm The frequency-domain prewhitening idea involved in
$T_n(\hat{\pmb{\gamma}})$ is similar in spirit to Hong's (1996),
where the prewhitening is done by fitting a long autoregression in
the time domain for each time series. However, Hong's test
statistic is valid only in the short memory case, while
$T_n(\hat{\pmb{\gamma}})$ is applicable to short memory,
anti-persistent and long memory time series. }
\end{remark}

 The theoretical results presented in this section state the
asymptotic null distributions of the infeasible test statistic
$T_n$ and two feasible ones $T_n(\hat{\pmb{\theta}})$ and
$T_n(\hat{{\pmb{\gamma}}})$. It would be desirable to study the
asymptotic distributions of these test statistics under either
fixed or local alternatives. At this point we are unable to obtain
any results in this direction due to some technical difficulties.
We conjecture that for fixed alternatives, the rate at which $T_n$
($T_n(\hat{\pmb{\theta}})$, $T_n(\hat{{\pmb{\gamma}}})$) diverges
is $n/\sqrt{B_n}$; compare Theorem 4 in Hong (1996) and Theorem 2
in Bouhaddioui and Roy (2006).  In addition, our proofs heavily
rely on the iid assumptions on the innovations $\{u_t\}$ and
$\{v_t\}$. It would be interesting to see to what extent these
assumptions can be relaxed and if the asymptotic distributions of
these test statistics are still valid for nonlinear models.


In the next section we examine the size and power properties of
our tests through Monte Carlo simulations.

\section{Simulation Studies}
\label{sec:3}

We consider the following two models:
\begin{eqnarray}
\label{eq:ar1}
(1-B)^{0.2}(1-0.5B)X_{1t}=u_t,~~~~(1-B)^{0.4}(1-0.5B)X_{2t}=v_t,
\end{eqnarray}
and
\begin{eqnarray}
\label{eq:ma1}
(1-B)^{-0.2}X_{1t}=(1+0.5B)u_t,~~~~(1-B)^{-0.4}X_{2t}=(1+0.5B)v_t,
\end{eqnarray}
 where $\{u_t\}$
($\{v_t\}$) are iid $N(0,1)$ random variables.
So in (\ref{eq:ar1}) ((\ref{eq:ma1})), $X_{1t}$ has
 moderate long memory  (antipersistence) and $X_{2t}$ possesses strong long memory (antipersistence).
 Two series lengths $(n=64,128)$ are investigated.

 The simulation of two FARIMA series from model (\ref{eq:ar1}) involves the following steps:
 \begin{enumerate}
 \item Generate two iid innovation sequences $\{u_t\}_{t=-4000,\cdots,n}$
and $\{v_t\}_{t=-4000,\cdots,n}$ with pre-specified correlations.
\item Simulate  two $\mbox{AR}(1)$ time series
$\tilde{X}_{1t}=0.5\tilde{X}_{1t-1}+u_t$ and
$\tilde{X}_{2t}=0.5\tilde{X}_{2t-1}+v_t$ recursively with the
first 1000 starting values subsequently discarded.

\item Generate two FARIMA series by applying a truncated
autoregression (corresponding to the fractional differencing) of
3000 lags [see Martin and Wilkins (1999)] to
$\{\tilde{X}_{1t}\}_{t=-3000,\cdots,n}$ and
$\{\tilde{X}_{2t}\}_{t=-3000,\cdots,n}$ respectively. Finally, we
only retain the last $n$ observations to ensure good
approximation.
\end{enumerate}

The simulation based on model (\ref{eq:ma1}) is similar except
that we generate two $\mbox{MA}(1)$ time series in the second
step. In the calculation of $T_n(\hat{\pmb{\theta}})$,  the
Whittle estimator $\hat{\pmb{\theta}}$ was obtained assuming
 an $\mbox{FAR}(1,d)$ model for each time series. Thus the assumption is incorrect if the data are generated from model (\ref{eq:ma1}). This will give us some idea of
 the consequences of model misspecification.  To calculate $T_n(\hat{\pmb{\gamma}})$,
 we
 also use Whittle estimator $\hat{\pmb{\gamma}}_k$ under the model $\mbox{FAR}(p_k,d)$ for $k=1,2$. We fix
 $p_1=p_2=3$ for $n=64$ and $p_1=p_2=5$ for $n=128$.


We consider the following three alternatives, which admit the same
forms as those in Hong (1996). Note that our alternatives are
closer to the null hypothesis than those used in his paper.
 \vskip 0.3cm
 Alternative 1: $\rho_{uv}(j)=0.05$ if $j=0$ and $0$ otherwise.
\vskip 0.3cm
 Alternative 2: $\rho_{uv}(j)=0.05$ if $j=0$,
$\sin(0.05\pi j)/(\pi j)$ if $1\le j\le 8$ and $0$ otherwise.
\vskip 0.3cm
 Alternative 3: $\rho_{uv}(j)=0.05$ if $j=3$ and $0$
otherwise.


In the simulation of Alternative 3, we generate $\{u_t,v_t\}_{t=-4000,n}$ by multiplying the square root of
its covariance matrix, which is a sparse matrix with dimension $(8002+2n)\times (8002+2n)$, with a vector of $(8002+2n)$
 random numbers independently generated from  standard normal distribution.  We use the following three kernel functions:
\begin{enumerate}
\item Bartlett (BAR) $K(x)=1-|x|$ if $|x|\le 1$ and $0$ otherwise;

\item Tukey (TUK) $K(x)=(1+\cos(\pi x))/2$, if $|x|\le 1$ and $0$
otherwise;

\item Parzen (PAR) $K(x)=1-6x^2+6|x|^3$, if $|x|\le 1/2$,
$2(1-|x|)^3$ if $1/2\le |x|\le 1$ and $0$ otherwise.
\end{enumerate}
For each kernel, we try three bandwidths: $B_n=[3n^{0.2}]$,
$[3n^{0.3} ]$ and $[3n^{0.4}]$, where $[a]$ stands for the integer
part of $a$. These rates lead to $B_n=6,10,15$ for $n=64$ and
$B_n=7,12,20$ for $n=128$.


Tables~\ref{tb:table1}-\ref{tb:table4} show the proportion of 5000
replications in which the null hypothesis was rejected at $5\%$
and $10\%$ nominal significance levels. Both tests are carried out
using upper tail critical values of the standard normal
distribution. Note that the estimated standard error of the relative rejection frequency
 is given by $\sqrt{\alpha(1-\alpha)/5000}$, where $\alpha$ is the observed relative rejection frequency.
 Under model (\ref{eq:ar1}), it is seen from Table 1 that
  both tests are oversized at $5\%$ and $10\%$ levels with less
size distortion for $n=128$. The larger bandwidth and the Tukey kernel correspond to
better size.
Under model (\ref{eq:ma1}), the size
distortion of $T_n(\hat{\pmb{\theta}})$ is noticeably larger than
that of $T_n(\hat{\pmb{\gamma}})$ uniformly in the kernel and
bandwidth. We attribute this phenomenon to the model
misspecification, which leads to inadequate prewhitening for
$T_n(\hat{\pmb{\theta}})$.
While for $T_n(\hat{\pmb{\gamma}})$, the $\mbox{FAR}(p_k,d)$,
$k=1,2$ models provide a decent approximation to the true data
generating process, so the size is much less distorted.

Tables~\ref{tb:table2}-\ref{tb:table4}  report size-adjusted power
of our test statistics under three alternatives. Under Alternative
1, the larger bandwidth corresponds to lower power. This is
expected since the test statistics assigning weights to a large
number of lags  are less powerful than those assigning  more
weights to short lags in detecting the simultaneous
cross-correlation of $\{u_t\}$ and $\{v_t\}$. Under model (\ref{eq:ar1}), the test statistic
$T_n(\hat{\pmb{\theta}})$ seems slightly more powerful than
$T_n(\hat{\pmb{\gamma}})$ for all three alternatives.
While for model (\ref{eq:ma1}), $T_n(\hat{\pmb{\gamma}})$ outperforms $T_n(\hat{\pmb{\theta}})$
in terms of power performance, which might be due to the model misspecification in forming $T_n(\hat{\pmb{\theta}})$.
Under Alternative 3, the power for $B_n=[3n^{0.2}]$ is
substantially lower than that for the other two bandwidths.
Compared to the other two kernels, the Parzen kernel performs very
poorly when $B_n=[3n^{0.2}]$. Since the correlation of $\{u_t\}$
and $\{v_t\}$ occurs only at lag $3$ under Alternative 3, any test
that gives less weight to the 3rd lag squared cross-correlation
tends to have less power. To illustrate this point, we note that
when $B_n=[3n^{0.2}]$, the (normalized) weights assigned to the
3rd lag for Bartlett, Tukey and Parzen kernels are $0.062$,
$0.056$ and $0.019$ for $n=64$, $0.069$, $0.071$ and $0.036$ for
$n=128$.
Since our
results only require finite fourth moment of the innovation
processes $\{u_t\}$ and $\{v_t\}$, we repeat the above simulations
with the innovations generated from a $t$ distribution with $5$
degrees of freedom. The performance of the tests (results not
shown) is very close to that of the tests when the innovations are
Gaussian.

Since $T_n(\hat{\pmb{\theta}})$ ($T_n(\hat{\pmb{\gamma}})$) is asymptotically equivalent to $G_n$ (see (\ref{eq:Gn})), which is a positive linear combination  of positive random
variables, its distribution is right-skewed and non-normal in the
small sample case. To correct for the size distortion observed in the simulation studies, we apply the power transformation method [Chen and Deo
(2004b)] and compare the performance of the modified statistic
with that of the original one (results not shown). It turns out that the power transformation method does not always lead to better size.  For Hong's (1996) test statistics, similar findings have  been reported in Chen and Deo (2004b), who attributed these to the fact that  the mean and variance of the test statistics
 are based on the asymptotic argument, thus could be inaccurate in small samples. In our setting, a detailed check of the sampling distribution of our test statistics shows that this is indeed the case for certain combinations of  the bandwidth and the kernel.  To obtain more accurate approximation to the mean and variance of our test statistics, one can follow the approach in Box and Pierce (1970) and Ling and Li (1997) to derive the asymptotic expression of $n\var(\hat{\pmb{\rho}})$, where $\hat{\pmb{\rho}}=(\hat{\rho}_{uv}(-B_n),\cdots,\hat{\rho}_{uv}(B_n))'$.   This is beyond the scope of this paper and will be an interesting topic for future research.

In summary, we observe reasonable size and power properties
for the proposed test statistics. The test statistics are oversized when the sample size $n=64$ and $128$, and the size distortion becomes less severe when the sample size is larger, e.x. $n=400$ (results not shown).  Under three alternatives
examined, the power is fairly high even for a moderate sample size
$n=128$ so long as the kernel and the bandwidth $B_n$ are appropriately chosen.
 Further research on the selection of the bandwidth parameter $B_n$  is needed and it seems to depend on the alternative under consideration.
Also no particular kernel was found to outperform the other two in
any situation considered.
Empirically, model misspecification was found to yield inferior performance in terms of  both size and power.
 A theoretical investigation of the
impact of the model misspecification on the asymptotic
distribution of $T_n(\hat{\pmb{\theta}})$ is certainly of
interest. We leave this topic for future work.
%

 \bigskip


 \baselineskip=17pt \centerline{REFERENCES}

\bigskip

\par\noindent\hangindent2.3em\hangafter 1
{Beran, J.} (1994) {\it Statistic for Long-Memory Processes},
Chapman \& Hall, New York, 1994.

\par\noindent\hangindent2.3em\hangafter 1
{Beran, J.}, { R. J. Bhansali} \& { D. Ocker} (1998) On unified
model selection for stationary and nonstationary short- and
long-memory autoregressive processes. {\it Biometrika} {85},
921-934.

\par\noindent\hangindent2.3em\hangafter 1
{Bhansali, R. J.}, { L. Giraitis} \& { P. S. Kokoszka} (2006)
Estimation of the memory parameter by fitting fractionally
differenced autoregressive models. {\it Journal of Multivariate
Analysis} {97}, 2101-2130.

\par\noindent\hangindent2.3em\hangafter 1
{Bouhaddioui, C.} \& { R. Roy} (2006) A generalized portmanteau
test for independence of two infinite-order vector autoregressive
series. {\it Journal of Time Series Analysis} {27}, 505-544.

\par\noindent\hangindent2.3em\hangafter 1
{Box, G.} \& {D. Pierce} (1970) Distribution of residual
autocorrelations in autoregressive-integrated moving average time
series models. {\it Journal of the American Statistical
Association}, 65, 1509-1526.

\par\noindent\hangindent2.3em\hangafter 1
{Brown, B.} (1971) Martingale central limit theorems. {\it Annals
of Mathematical Statistics} 42, 59-66.


\par\noindent\hangindent2.3em\hangafter 1
{Chen, W.} \& { R. S. Deo} (2004a) A generalized portmanteau
goodness-of-fit test for time series models. {\it Econometric
Theory} {\bf 20} 382-416.

\par\noindent\hangindent2.3em\hangafter 1
{Chen, W.} \& { R. S. Deo} (2004b) Power transformations  to induce normality
 and their applications.  {\it Journal of Royal Statistical Society, B.} {\bf 66}
 117-130.

\par\noindent\hangindent2.3em\hangafter 1
{Dahlhaus, R.} (1989) Efficient parameter estimation for
self-similar processes. {\it Annals of Statistics} 17, 1749-1766.


\par\noindent\hangindent2.3em\hangafter 1
{Doukhan, P.}, { G. Oppenheim} \& { M. S. Taqqu} (2003) {\it
Theory and Applications of Long-range Dependence}, Birkha\"user,
Boston.

\par\noindent\hangindent2.3em\hangafter 1
{Duchesne P.} \& { R. Roy} (2003) Robust tests for independence
of two time series. {\it Statistica Sinica} {13}, 827-852.


\par\noindent\hangindent2.3em\hangafter 1
{Eichler, M.} (2007) A frequency-domain based test for
non-correlation between stationary time series. {\it Metrika}, 65, 133-157.

\par\noindent\hangindent2.3em\hangafter 1
{El Himdi, K.} \& { R. Roy} (1997) Tests for noncorrelation of
two multivariate ARMA time series.  {\it Canadian Journal of
Statistics} {25}, 233-256.

\par\noindent\hangindent2.3em\hangafter 1
{Fox, R.} \& {M. S. Taqqu} (1986) Large-sample properties of
parameter estimates for strongly dependent stationary Gaussian
time series. {\it Annals of Statistics} 14, 517-532.

\par\noindent\hangindent2.3em\hangafter 1
{Giraitis, L.} \& {D. Surgailis} (1990) A central limit theorem
for quadratic forms in strongly dependent random variables and its
application to asymptotic normality of Whittle's estimate. {\it
Probability Theory and Related Fields} 86, 87-104.


\par\noindent\hangindent2.3em\hangafter 1
{Hallin, M.} \& {A. Saidi} (2005) Testing non-correlation and
non-causality between multivariate ARMA time series.  {\it Journal
of Time Series Analysis} 26, 83-106.

\par\noindent\hangindent2.3em\hangafter 1
{Hallin, M.} and {A. Saidi} (2007) Optimal tests of non-correlation between
multivariate time series. {\it Journal of the American Statistical Association}, 102, 938-952.

\par\noindent\hangindent2.3em\hangafter 1
{Hannan, E. J.} (1973) The asymptotic theory of linear time
series models. {\it Journal of Applied Probability} 10, 130-145.

\par\noindent\hangindent2.3em\hangafter 1
{Haugh, L. D.} (1976) Checking the independence of two
covariance-stationary time series: a univariate residual
cross-correlation approach. {\it Journal of the American
Statistical Association} 71, 378-385.

\par\noindent\hangindent2.3em\hangafter 1
{Hong, Y.} (1996) Testing for independence between two covariance
stationary time series. {\it Biometrika} {83}, 615-625.

\par\noindent\hangindent2.3em\hangafter 1
{Ling, S.} and {W. K. Li} (1997) On fractionally integrated
autoregressive moving-average time series models with conditional
heteroscedasticity. {\it Journal of the American Statistical
Association}, 92, 1184-1194.





\par\noindent\hangindent2.3em\hangafter 1
{Martin, V. L.} \& { N. P. Wilkins} (1999). Indirect estimation
of ARFIMA and VARFIMA models. {\it Journal of Econometrics} {93},
149-175.

\par\noindent\hangindent2.3em\hangafter 1
{Pham, D.}, { R. Roy} \& { L. C\'{e}dras} (2003) Tests for
non-correlation of two cointegrated ARMA time series. {\it Journal
of Time Series Analysis} 24, 553-577.

\par\noindent\hangindent2.3em\hangafter 1
{Priestley, M. B.} (1981) {\it Spectral Analysis and Time
Series.} Volume 1, Academic, New York.

\par\noindent\hangindent2.3em\hangafter 1
{Robinson, P. M.} (1995) Gaussian semiparametric estimation of
long range dependence.  {\it Annals of Statistics} 23, 1630-1661.

\par\noindent\hangindent2.3em\hangafter 1
{Robinson, P. M.} (2003) {\it Time Series with Long Memory},
Oxford University Press.


\par\noindent\hangindent2.3em\hangafter 1
{Teyssi$\grave{e}$re, G.} \& {Kirman, A.} (2005) {\it Long Memory
in Economics}, Springer Verlag.

\par\noindent\hangindent2.3em\hangafter 1
{Velasco, C.} \& { P. M. Robinson} (2000) Whittle pseudo-maximum
likelihood estimation for nonstationary time series. {\it Journal
of the American Statistical Association} 95, 1229-1243.

\par\noindent\hangindent2.3em\hangafter 1
{Wu, W. B.} \& { X. Shao} (2007). A limit theorem for quadratic
forms and its applications. {\it Econometric Theory}, 23, 930-951.

\bigskip

\section{Appendix}
In the appendix, the constant $C$ is generic and it may vary from
line to line.  The following decompositions (\ref{eq:bartlettX})
and (\ref{eq:bartlettY}) will be frequently used in the proof.
Note that
\begin{eqnarray*}
w_{X_1}(\lambda)=\frac{1}{\sqrt{2\pi
n}}\sum_{t=1}^{n}\left(\sum_{j=0}^{\infty}a_j u_{t-j}\right)
e^{it\lambda}=A(\lambda)w_{u}(\lambda)+\frac{1}{\sqrt{2\pi
n}}\sum_{j=0}^{\infty} a_j U_{j,n}(\lambda) e^{ij\lambda},
\end{eqnarray*}
where $U_{j,n}(\lambda)=\sum_{t=1-j}^{n-j} u_t
e^{it\lambda}-\sum_{t=1}^{n} u_t e^{it\lambda}$. Since
$f_{X_1X_1}(\lambda)=(2\pi)^{-1}|A(\lambda)|^2$, we have
\begin{eqnarray}
\label{eq:bartlettX}
\frac{w_{X_1}(\lambda)}{\sqrt{f_{X_1X_1}(\lambda)}}=\frac{A(\lambda)
w_{u}(\lambda)}{|A(\lambda)|/\sqrt{2\pi}}+R_u(\lambda),
~\mbox{where}~R_u(\lambda)=\frac{\sum_{j=0}^{\infty} a_j
U_{j,n}(\lambda) e^{ij\lambda}}{\sqrt{n}|A(\lambda)|}.
\end{eqnarray}
Similarly, we have
\begin{eqnarray}
\label{eq:bartlettY}
\frac{w_{X_2}(\lambda)}{\sqrt{f_{X_2X_2}(\lambda)}}=\frac{B(\lambda)
w_{v}(\lambda)}{|B(\lambda)|/\sqrt{2\pi}}+R_v(\lambda),~\mbox{where}~
R_v(\lambda)=\frac{\sum_{j=0}^{\infty} b_j V_{j,n}(\lambda)
e^{ij\lambda}}{\sqrt{n}|B(\lambda)|}
\end{eqnarray}
and $V_{j,n}(\lambda)=\sum_{t=1-j}^{n-j} v_t
e^{it\lambda}-\sum_{t=1}^{n} v_t e^{it\lambda}$.

For the convenience of notation, write $A_j=A(\lambda_j)$,
$B_j=B(\lambda_j)$, $I_{12j}=I_{X_1X_2}(\lambda_j)$,
$w_{uj}=w_{u}(\lambda_j)$, $w_{vj}=w_{v}(\lambda_j)$,
$R_{uj}=R_{u}(\lambda_j)$ and $R_{vj}=R_{v}(\lambda_j)$.
Let $g_j=A_j|A_j|^{-1}$ and $h_j=B_j|B_j|^{-1}$. Further let
$D_k(\lambda)=\sum_{t=1}^{k} e^{it\lambda}$. Denote by $a\vee
b=\max(a,b)$ and $a\wedge b=\min(a,b)$. Define
\[\hat{f}_{uv}(\lambda)=\frac{2\pi}{n}\sum_{j=1}^{n-1}
W(\lambda-\lambda_j)g_jw_{u}(\lambda_j)\overline{h_j
w_{v}(\lambda_j)},\]
\[\hat{f}_{uu}(\lambda)=\frac{2\pi}{n}\sum_{j=1}^{n-1}
W(\lambda-\lambda_j)I_{uu}(\lambda_j)~\mbox{and
}~\hat{f}_{vv}(\lambda)=\frac{2\pi}{n}\sum_{j=1}^{n-1}
W(\lambda-\lambda_j)I_{vv}(\lambda_j).\]

 \noindent Proof of
Theorem~\ref{th:mainresult1}: The conclusion follows from the
following three assertions:
\begin{eqnarray*}
\label{eq:neg10}
\sum_{l=0}^{n-1}\left(|\hat{f}_{X_1X_2}(\lambda_l)|^2-4\pi^2|\hat{f}_{uv}(\lambda_l)|^2\right)=o_{p}(\sqrt{B_n}),
\end{eqnarray*}
\begin{eqnarray*}
\label{eq:neg2} \frac{1}{n}\sum_{l=0}^{n-1}
\hat{f}_{X_kX_k}(\lambda_l)-1=o_{p}(B_n^{-1/2}),~k=1,2
\end{eqnarray*}
\begin{eqnarray*}
\label{eq:clt}
~\mbox{and}~~~~~~~~\frac{4\pi^2\sum_{l=0}^{n-1}|\hat{f}_{uv}(\lambda_l)|^2-B_ns(K)}{\sqrt{2B_n
d(K)}}\rightarrow_{D} N(0,1),
\end{eqnarray*}
which have been established in
Lemmas~\ref{lem:numerator},~\ref{lem:denominatorneg} and
\ref{lem:clt} respectively. \qed

\begin{lemma}
\label{lem:numerator} Under the assumptions of
Theorem~\ref{th:mainresult1}, it follows that
\begin{eqnarray}
\label{eq:neg10}
\sum_{l=0}^{n-1}\left(|\hat{f}_{X_1X_2}(\lambda_l)|^2-4\pi^2|\hat{f}_{uv}(\lambda_l)|^2\right)=o_{p}(\sqrt{B_n}).
\end{eqnarray}
\end{lemma}
\noindent Proof of Lemma~\ref{lem:numerator}: According to
(\ref{eq:bartlettX}) and (\ref{eq:bartlettY}), we get
\begin{eqnarray*}
\hat{f}_{X_1X_2}(\lambda_l)-2\pi
\hat{f}_{uv}(\lambda_l)=\frac{2\pi}{n}\sum_{j=1}^{n-1}W(\lambda_{l-j})[\sqrt{2\pi}
g_j
w_{uj}\overline{R_{vj}}+\sqrt{2\pi}\overline{h_jw_{vj}}R_{uj}+R_{uj}\overline{R_{vj}}].
\end{eqnarray*}
Then $\mbox{LHS (left hand side)
of}~(\ref{eq:neg10})=G_{1n}+G_{2n}+\overline{G_{2n}}$, where
\begin{eqnarray*}
G_{1n}&=&\frac{4\pi^2}{n^2}
\sum_{l=0}^{n-1}\left|\sum_{j=1}^{n-1}W(\lambda_{l-j})[\sqrt{2\pi}g_j
w_{uj}\overline{R_{vj}}+\sqrt{2\pi}\overline{h_jw_{vj}}R_{uj}+R_{uj}\overline{R_{vj}}]\right|^2,\\
G_{2n}&=&\frac{4\pi^2}{n}\sum_{l=0}^{n-1}\overline{\hat{f}_{uv}(\lambda_l)}\sum_{j=1}^{n-1}W(\lambda_{l-j})[\sqrt{2\pi}
g_jw_{uj}\overline{R_{vj}}+\sqrt{2\pi}\overline{h_jw_{vj}}R_{uj}+R_{uj}\overline{R_{vj}}].
\end{eqnarray*}
To show $G_{1n}=o_{p}(\sqrt{B_n})$, we note that $G_{1n}$ is
smaller than $4(G_{11n}+G_{12n}+G_{13n})$, where
\begin{eqnarray*}
G_{11n}&=&\frac{8\pi^3}{n^2}\sum_{l=0}^{n-1}\sum_{j,j'=1}^{n-1}W(\lambda_{l-j})W(\lambda_{l-j'})g_jw_{uj}\overline{R_{vj}} \overline{g_{j'}}\overline{w_{uj'}}R_{vj'},\\
G_{12n}&=&\frac{8\pi^3}{n^2}\sum_{l=0}^{n-1}\sum_{j,j'=1}^{n-1}W(\lambda_{l-j})W(\lambda_{l-j'})\overline{h_j}\overline{w_{vj}}R_{uj}h_{j'}w_{vj'}\overline{R_{uj'}},\\
G_{13n}&=&\frac{4\pi^2}{n^2}\sum_{l=0}^{n-1}\sum_{j,j'=1}^{n-1}W(\lambda_{l-j})W(\lambda_{l-j'})R_{uj}\overline{R_{vj}}
\overline{R_{uj'}} R_{vj'}.
\end{eqnarray*}
We shall only prove $G_{11n}=o_p(\sqrt{B_n})$, since the treatment
for $G_{12n}$ and $G_{13n}$ is similar. Let
$\Lambda_n(\lambda_s)=\sum_{h=-B_n}^{B_n}K^2(hb_n)e^{ih\lambda_s}$,
where $b_n=B_n^{-1}$. In view of the fact that
\begin{eqnarray}
\label{eq:WLambda}
\sum_{l=0}^{n-1}W(\lambda_{l-j})W(\lambda_{l-j'})=(4\pi^2)^{-1} n
\Lambda_n(\lambda_{j-j'}),
\end{eqnarray}
we have
\[G_{11n}= 2\pi
n^{-1}\sum_{j,j'=1}^{n-1}g_j\overline{g_{j'}}w_{uj}\overline{w_{uj'}}\overline{R_{vj}}R_{vj'}\Lambda_n(\lambda_{j-j'}).\]
By Lemma~\ref{lem:wwww}, we get
\begin{eqnarray*}
\E|G_{11n}|^2&=&\frac{4\pi^2}{n^{2}}\sum_{j_1,j_1'=1}^{n-1}\sum_{j_2,j_2'=1}^{n-1}
\E[w_{uj_1}\overline{w_{uj_1'}}\overline{w_{uj_2}}{w_{uj_2'}}]\E[\overline{R_{vj_1}}
R_{v j_1'}
R_{vj_2}\overline{R_{vj_2'}}]\\
&&g_{j_1}\overline{g_{j_1'}}\overline{g_{j_2}}{g_{j_2'}}\Lambda_n(\lambda_{j_1-j_1'})\Lambda_n(\lambda_{j_2-j_2'})\\
&\le&\frac{1}{n^2}\sum_{j_1,j_2=1}^{n-1}\Lambda_n^2(0)|\E[\overline{R_{vj_1}}
R_{v j_1}
R_{vj_2}\overline{R_{vj_2}}]|+\frac{1}{n^2}\sum_{j_1,j_1'=1}^{n-1}\Lambda_n^2(\lambda_{j_1-j_1'})\\
&&\times\{|\E[\overline{R_{vj_1}} R_{v j_1'}
{R_{vj_1}}\overline{R_{vj_1'}}]|+|\E[\overline{R_{vj_1}}
R_{v j_1'} R_{vj_1}\overline{R_{vj_1'}}]|\}\\
&&+\frac{|c_4(u)-3|}{n^3}\sum_{j_1,j_1',j_2=1}^{n-1}\Lambda_n^2(\lambda_{j_1-j_1'})|\E[\overline{R_{vj_1}}
R_{v j_1'} R_{vj_2}\overline{R_{v(j_1+j_2-j_1')}}]|.
\end{eqnarray*}
Let ${\bf 1}(\cdot)$ denotes the indicator function.  The
following fact
\begin{eqnarray}
\label{eq:Lambdan} |\Lambda_n(\lambda_s)|={\bf
1}(s=0~\mbox{mod}~n)O(B_n)+O(|\lambda_s|^{-1}){\bf
1}(s\not=0~\mbox{mod}~n),
\end{eqnarray}
which was stated in Equation A.28 of Chen and Deo (2004), and
Lemma~\ref{lem:bound1} result in  $\E|G_{1n}|^2=O(\log^2 n
B_n^2/n^2)+O(B_n^2/n^{3/2})=o({B_n})$. Hence
$G_{1n}=o_p(\sqrt{B_n})$.

 To
show $G_{2n}=o_{p}(\sqrt{B_n})$, we first show that
\begin{eqnarray*}
\label{eq:G11n}
G_{21n}:=\frac{4\pi^2\sqrt{2\pi}}{n}\sum_{l=0}^{n-1}\overline{\hat{f}_{uv}(\lambda_l)}\sum_{j=1}^{n-1}W(\lambda_{l-j})
g_j w_{uj}\overline{R_{vj}}=o_{p}(\sqrt{B_n}).
\end{eqnarray*}
  Using (\ref{eq:WLambda})
again, we get
\begin{eqnarray*}
G_{21n}&=&\frac{8\pi^3\sqrt{2\pi}}{n^2}\sum_{l=0}^{n-1}\sum_{j,j'=1}^{n-1}W(\lambda_{l-j})W(\lambda_{l-j'})g_j\overline{g_{j'}}h_{j'}w_{uj}\overline{w_{uj'}}w_{vj'}\overline{R_{vj}}\\
       &=&\frac{2\pi\sqrt{2\pi}}{n}\sum_{j,j'=1}^{n-1}\Lambda_n(\lambda_{j'-j})g_j\overline{g_{j'}}h_{j'} w_{uj}\overline{w_{uj'}}w_{vj'}\overline{R_{vj}}.
\end{eqnarray*}
So
\begin{eqnarray*}
\E|G_{21n}|^2&=&\frac{8\pi^3}{n^2}\sum_{j_1,j_1'=1}^{n-1}\sum_{j_2,j_2'=1}^{n-1}\E(w_{uj_1}\overline{w_{uj_1'}}\overline{w_{uj_2}}w_{uj_2'})\E(w_{vj_1'}\overline{R_{vj_1}}
\overline{w_{vj_2'}}R_{vj_2})\\
&&\times
\Lambda_n(\lambda_{j_1'-j_1})g_{j_1}\overline{g_{j_1'}}h_{j_1'}\Lambda_n(\lambda_{j_2'-j_2})\overline{g_{j_2}}{g_{j_2'}}\overline{h_{j_2'}}.
\end{eqnarray*}
Similar to $\E|G_{11n}|^2$ above,  $\E|G_{21n}|^2$ can be bounded
by a sum of four terms in view of Lemma~\ref{lem:wwww}.  For
example, the first term corresponds to the case $j_1=j_1'$,
$j_2=j_2'$, which is
\[O(n^{-2})\sum_{j_1,j_2=1}^{n-1}\Lambda_n^2(0)|\E(w_{vj_1}\overline{R_{vj_1}}
\overline{w_{vj_2}}R_{vj_2})|=O(B_n^2n^{-2})\sum_{j_1,j_2=1}^{[n/2]}(j_1j_2)^{-1/2}=o(B_n)\]
 by Lemma~\ref{lem:bound1}.
 The bound for the other three terms can be
established using Lemma~\ref{lem:bound1} in a similar fashion.
 Hence $G_{21n}=o_{p}(\sqrt{B_n})$. Following the same argument,
 the other two terms in $G_{2n}$ can be shown to be
 $o_p(\sqrt{B_n})$. This completes the proof.

\qed

\begin{lemma}
\label{lem:denominatorneg} Under the assumptions of
Theorem~\ref{th:mainresult1}, $n^{-1}\sum_{l=0}^{n-1}
\hat{f}_{X_kX_k}(\lambda_l)-1=o_{p}(B_n^{-1/2})$, $k=1,2$.
\end{lemma}

\noindent Proof of Lemma~\ref{lem:denominatorneg}: We only deal
with the case $k=1$. Observing that
\begin{eqnarray*}
\label{eq:Wlambda} \sum_{l=0}^{n-1}W(\lambda_{l-j})=(2\pi)^{-1}n~
\mbox{for}~ j=1,2,\cdots,n-1,
\end{eqnarray*}
we derive
\begin{eqnarray*}
&&\hspace{-0.5cm}\frac{1}{n}\sum_{l=0}^{n-1}
\hat{f}_{X_1X_1}(\lambda_l)=\frac{2\pi}{n^2}\sum_{l=0}^{n-1}\sum_{j=1}^{n-1}\frac{W(\lambda_{l-j})I_{X_1X_1}(\lambda_j)}{f_{X_1X_1}(\lambda_j)}=
\frac{1}{n}\sum_{j=1}^{n-1}[2\pi
I_{uu}(\lambda_j)+|R_{uj}|^2\\
&&\hspace{0.5cm}+\sqrt{2\pi}g_jw_{uj}\overline{R}_{uj}+\sqrt{2\pi}\overline{g_j}\overline{w_{uj}}{R}_{uj}]=:J_1+J_2+J_3+\overline{J_3}.
\end{eqnarray*}
 Let
$\bar{I}_{uu}(\lambda)=(2\pi
n)^{-1}|\sum_{t=1}^{n}(u_t-\bar{u}_n)e^{it\lambda}|^2$ be the
periodogram of the centered innovations $u_t-\bar{u}_n$ with
$\bar{u}_n=n^{-1}\sum_{t=1}^{n}u_t$. Since
$I_{uu}(\lambda_j)=\bar{I}_{uu}(\lambda_j)$ at $j=1,2,\cdots,n-1$
and $2\pi
n^{-1}\sum_{j=0}^{n-1}\bar{I}_{uu}(\lambda_j)=n^{-1}\sum_{t=1}^n(u_t-\bar{u}_n)^2$,
we have
$J_1=n^{-1}\sum_{t=1}^n(u_t-\bar{u}_n)^2=1+O_{p}(n^{-1/2})$ by the
central limit theorem.

By Lemma~\ref{lem:bound1}, we get
\begin{eqnarray*}
\var(J_3)&=&\frac{2\pi}{n^2}\sum_{j,j'=1}^{n-1}g_j\overline{g_{j'}}\cov(w_{uj}\overline{R_{uj}},\overline{w_{uj'}}{R_{uj'}})\\
         &=&\frac{2\pi}{n^2}\sum_{j,j'=1}^{n-1}g_j\overline{g_{j'}}[\cum(w_{uj},\overline{R_{uj}},\overline{w_{uj'}},{R_{uj'}})+\cov(w_{uj},\overline{w_{uj'}})\cov(\overline{R_{uj}},R_{uj'})\\
         &&+\cov(w_{uj},{R_{uj'}})\cov(\overline{w_{uj'}},\overline{R_{uj}})]=O(n^{-2})\sum_{j,j'=1}^{[n/2]}(jj')^{-1/2}=O(n^{-1}).
\end{eqnarray*}
 A similar argument yields $\var(J_2)=O(n^{-2}\log^2 n)$.
The conclusion follows.

\qed

\begin{lemma}
\label{lem:clt} Under the assumptions of
Theorem~\ref{th:mainresult1}, we have that
\[\frac{4\pi^2\sum_{l=0}^{n-1}|\hat{f}_{uv}(\lambda_l)|^2-B_ns(K)}{\sqrt{2B_n
d(K)}}\rightarrow_{D} N(0,1).\]
\end{lemma}

\noindent Proof of Lemma~\ref{lem:clt}:  Let
$G_{j,j'}:=g_j\overline{h_j}\overline{g_{j'}} h_{j'}$ and
$H_n:=4\pi^2\sum_{l=0}^{n-1}|\hat{f}_{uv}(\lambda_l)|^2$. Then by
(\ref{eq:WLambda}),
\begin{eqnarray}
\label{eq:Hn}
H_n&=&\frac{16\pi^4}{n^2}\sum_{l=0}^{n-1}\sum_{j,j'=1}^{n-1}W(\lambda_{l-j})W(\lambda_{l-j'})g_jw_{uj}\overline{h_jw_{vj}}
\cdot \overline{g_{j'}w_{uj'}}h_{j'}w_{vj'}\nonumber\\
&=&\frac{4\pi^2}{n}\sum_{j,j'=1}^{n-1}G_{j,j'}\Lambda_n(\lambda_{j-j'})w_{uj}\overline{w_{vj}}\overline{w_{uj'}}w_{vj'}\\
   &=&\frac{1}{n^3}\sum_{t_1=1}^{n}\sum_{t_2=1}^{n}\sum_{t_3=1}^{n}\sum_{t_4=1}^{n}u_{t_1}u_{t_2}v_{t_3}v_{t_4}
   \sum_{j,j'=1}^{n-1}G_{j,j'}\Lambda_n(\lambda_{j-j'}) e^{i (t_1-t_3)\lambda_j}
   e^{-i(t_2-t_4)\lambda_{j'}}.\nonumber
\end{eqnarray}
Denote by
$g_n(k_1,k_2)=\sum_{j,j'=1}^{n-1}G_{j,j'}\Lambda_n(\lambda_{j-j'})
e^{i k_1\lambda_j} e^{-i k_2\lambda_{j'}}$ and $g_n(k)=g_n(k,k)$.
Note that $g_n(k)$ is real. Write
\begin{eqnarray*}
&&\hspace{-0.6cm}H_{n}=\frac{1}{n^3}\sum_{k_1,k_2=1-n}^{n-1}\sum_{t_1=(1+k_1)\vee
1}^{(n+k_1)\wedge n}\sum_{t_2=(1+k_2)\vee 1}^{(n+k_2)\wedge
n}u_{t_1}
u_{t_2} v_{t_1-k_1}v_{t_2-k_2}g_n(k_1,k_2)\\
&&\hspace{-0.5cm}=H_{0n}+\frac{1}{n^3}\sum_{k=1-n}^{n-1}\sum_{t_1=(1+k)\vee
1}^{(n+k)\wedge n}\sum_{t_2=(1+k)\vee 1}^{(n+k)\wedge n}u_{t_1}
u_{t_2} v_{t_1-k}v_{t_2-k}g_n(k)=H_{0n}+H_{1n}+H_{2n},
\end{eqnarray*}
where
\begin{eqnarray*}
H_{0n}&=&\frac{1}{n^3}\sum_{k_1\not
=k_2=1-n}^{n-1}\sum_{t_1=(1+k_1)\vee 1}^{(n+k_1)\wedge
n}\sum_{t_2=(1+k_2)\vee 1}^{(n+k_2)\wedge n}u_{t_1}
u_{t_2} v_{t_1-k_1}v_{t_2-k_2}g_n(k_1,k_2),\\
H_{1n}&=&\frac{1}{n^3}\sum_{k=0}^{n-1}\sum_{t_1,t_2=1+k}^{n}u_{t_1}
u_{t_2} v_{t_1-k}v_{t_2-k}g_n(k),\\
H_{2n}&=&
\frac{1}{n^3}\sum_{k=1}^{n-1}\sum_{t_1,t_2=k+1}^{n}v_{t_1}v_{t_2}u_{t_1-k}u_{t_2-k}g_n(-k).
\end{eqnarray*}
We shall first prove $H_{0n}=o_p(\sqrt{B_n})$. Note that
\begin{eqnarray*}
\var(H_{0n})&=&n^{-6}\sum_{k_1\not =k_2=1-n}^{n-1} \sum_{k_1'\not
=k_2'=1-n}^{n-1}\sum_{t_1=(1+k_1)\vee 1}^{(n+k_1)\wedge
n}\sum_{t_2=(1+k_2)\vee 1}^{(n+k_2)\wedge n}
\sum_{t_1'=(1+k_1')\vee 1}^{(n+k_1')\wedge n}\sum_{t_2'=(1+k_2')\vee 1}^{(n+k_2')\wedge n}\\
&& \E[u_{t_1} u_{t_2} u_{t_1'} u_{t_2'}] \E[v_{t_1-k_1}v_{t_2-k_2}
v_{t_1'-k_1'}v_{t_2'-k_2'}]g_n(k_1,k_2)\overline{g_n(k_1',k_2')}.
\end{eqnarray*}
By Lemma~\ref{lem:wwww}, there are at least four restrictions on
$t_1,t_2,t_1',t_2',k_1,k_2,k_1',k_2'$ for each non-vanishing term
in $\var(H_{0n})$. For example, one of such terms, denoted by
$J_{1n}$, corresponds to $t_1=t_2,t_1'=t_2'$, $k_1=k_1',k_2=k_2'$,
i.e.
\begin{eqnarray*}
J_{1n}&=&O(n^{-6})\sum_{k_1\not
=k_2=1-n}^{n-1}\sum_{t_1,t_1'=(1+k_1)\vee (1+k_2)\vee
1}^{(n+k_1)\wedge (n+k_2)\wedge n}|g_n(k_1,k_2)|^2\\
&=&O(n^{-4})\sum_{k_1,k_2=1-n}^{n-1} |g_n(k_1,k_2)|^2\\
&=&O(n^{-4})\sum_{j_1,j_2,j_1',j_2'=1}^{n-1}|\Lambda_n(\lambda_{j_1-j_1'})\Lambda_n(\lambda_{j_2-j_2'})|\left|\sum_{k_1,k_2=1-n}^{n-1}e^{i(k_1\lambda_{j_1-j_2}-k_2\lambda_{j_1'-j_2'})}\right|.
\end{eqnarray*}
In view of (\ref{eq:Lambdan}) and the fact that
$\sum_{k=1-n}^{n-1}e^{ik\lambda_{s}}={\bf
1}(s=0~\mbox{mod}~n)(2n-1)-{\bf 1}(s\not=0~\mbox{mod}~n)$ for
$s\in\Z$, we have $J_{1n}=O(\log^2 n)=o(B_n)$ under
Assumption~\ref{as:bandwidth}. Other terms in $\var(H_{0n})$ can
be bounded by $O(\log^2 n)$ in a similar way, so
$\var(H_{0n})=o(B_n)$ and $H_{0n}=o_p(\sqrt{B_n})$.

 Write $H_{1n}=H_{11n}+H_{12n}$ and $H_{2n}=H_{21n}+H_{22n}$, where
 \begin{eqnarray*}
H_{11n}&=&\frac{2}{n^{3}}\sum_{k=0}^{n-1}\sum_{t_1=k+2}^{n}\sum_{t_2=k+1}^{t_1-1}u_{t_1}
u_{t_2}
v_{t_1-k}v_{t_2-k}g_n(k),\\
H_{21n}&=&\frac{2}{n^{3}} \sum_{k=1}^{n-1}
\sum_{t_1=k+2}^{n}\sum_{t_2=k+1}^{t_1-1} v_{t_1} v_{t_2}
u_{t_1-k}u_{t_2-k}g_n(-k),\\
H_{12n}&=&\frac{1}{n^{3}}\sum_{k=0}^{n-1}\sum_{t=k+1}^{n} u_{t}^2
v_{t-k}^2g_n(k),~H_{22n}=\frac{1}{n^{3}}\sum_{k=1}^{n-1}\sum_{t=k+1}^{n}v_{t}^2
u_{t-k}^2g_n(-k).
\end{eqnarray*}
We shall show that $\var(H_{12n})=o(B_n)$ and
$\var(H_{22n})=o(B_n)$.
 Note that
\begin{eqnarray*}
\var(H_{22n})&=&\frac{1}{n^6}\sum_{k_1,k_2=1}^{n-1}\sum_{t_1=k_1+1}^n\sum_{t_2=k_2+1}^{n}g_n(-k_1)g_n(-k_2)\cov(v_{t_1}^2u_{t_1-k_1}^2,v_{t_2}^2
u_{t_2-k_2}^2).
\end{eqnarray*}
Since $\cov(v_{t_1}^2u_{t_1-k_1}^2,v_{t_2}^2
u_{t_2-k_2}^2)=\E(v_{t_1}^2v_{t_2}^2)\E(u_{t_1-k_1}^2u_{t_2-k_2}^2)-1=[c_4(v)+2]{\bf
1}(t_1=t_2,k_1\not=k_2)+[c_4(u)+2]{\bf
1}(t_1-k_1=t_2-k_2,t_1\not=t_2)+[(c_4(v)+3)(c_4(u)+3)-1]{\bf
1}(t_1=t_2,k_1=k_2)$, we can write $\var(H_{22n})$ into  a sum of
three terms, which correspond to $t_1=t_2$, $t_1-k_1=t_2-k_2$
 and $t_1=t_2,k_1=k_2$. When $t_1=t_2$, the corresponding term is
\begin{eqnarray*}
J_{2n}&=&O(n^{-6})\sum_{k_1,k_2=1}^{n-1}\sum_{t_1=(k_1+1)\vee
(k_2+1)}^n
g_n(-k_1)g_n(-k_2)=O(n^{-6})\sum_{j_1,j_1',j_2,j_2'=1}^{n-1}\\
&&|\Lambda_n(\lambda_{j_1-j_1'})||\Lambda_n(\lambda_{j_2-j_2'})|
\left|\sum_{k_1,k_2=1}^{n-1}\sum_{t_1=(k_1+1)\vee (k_2+1)}^n
e^{-i(k_1\lambda_{j_1-j_1'}+k_2\lambda_{j_2-j_2'})}\right|\\
&=&O(B_n^2n^{-6})\sum_{j_1,j_1',j_2,j_2'=1}^{n-1}\left|\sum_{t_1=2}^{n}D_{t_1-1}(\lambda_{j_1-j_1'})
D_{t_1-1}(\lambda_{j_2-j_2'})\right|\\
&=&O(B_n^2n^{-4})\sum_{m_1,m_2=2-n}^{n-2}\left|\sum_{t_1=2}^{n}D_{t_1-1}(\lambda_{m_1})
D_{t_1-1}(\lambda_{m_2})\right|.
\end{eqnarray*}
Since $D_k(\lambda_s)=O(|\lambda_s|^{-1}){\bf 1}(s\not=0)+k{\bf
1}(s=0)$ for $s=2-n,3-n,\cdots,n-3,n-2$, under
Assumption~\ref{as:bandwidth}, we have
\begin{eqnarray*}
|J_{2n}|&=&O(B_n^2/n^{3})\sum_{m_1,m_2=1}^{n-2}(\lambda_{m_1}\lambda_{m_2})^{-1}+O(B_n^2
/n^2)\sum_{m_1=1}^{n-2}\lambda_{m_1}^{-1}\\
&&+O(B_n^2/n^{4})\sum_{t_1=2}^{n}(t_1-1)^2=O(B_n^2\log^2
n/n)=o(B_n).
\end{eqnarray*}
Similarly, we can show that the other two terms corresponding to
$t_1-k_1=t_2-k_2$ and $t_1=t_2,k_1=k_2$
 are of order $o(B_n)$. So $H_{22n}=o_{p}(\sqrt{B_n})$.
 Regarding $H_{12n}$, we note that $H_{12n}=n^{-3}\sum_{t=1}^{n}u_t^2
v_{t}^2g_n(0)+n^{-3}\sum_{k=1}^{n-1}\sum_{t=k+1}^{n} u_t^2
v_{t-k}^2g_n(k)$, where the variance of the latter term is of
order $o(B_n)$ by the same argument as above and the former term
is easily shown to be $o_p(\sqrt{B_n})$. Therefore
$H_n=H_{11n}+H_{21n}+o_p(\sqrt{B_n})$.

Since $\E(H_n)=B_ns(K)+O(1)$ and $\var(H_n)=2B_n d(K)(1+o(1))$, as
proved in Lemma~\ref{lem:variance},  our conclusion holds if we
can show
\begin{eqnarray}
\label{eq:1121}
\frac{H_{11n}+H_{21n}}{\sqrt{\sigma_n^2}}\rightarrow_{D}
N(0,1),~\mbox{where}~\sigma_n^2=2B_n d(K).
\end{eqnarray}
 Note that
$\sum_{k=1}^{n-1}\sum_{t=k+2}^{n}\sum_{s=k+1}^{t-1}=\sum_{t=3}^n\sum_{s=2}^{t-1}\sum_{k=1}^{s-1}$
and
$\sum_{k=0}^{n-1}\sum_{t=k+2}^{n}\sum_{s=k+1}^{t-1}=\sum_{t=2}^n$
$\sum_{s=1}^{t-1}\sum_{k=0}^{s-1}$. We can write
$H_{11n}=n^{-1}\sum_{t=2}^n W_{1nt}$ and
$H_{21n}=n^{-1}\sum_{t=3}^{n}W_{2nt}$, where
\[W_{1nt}=\frac{2}{n^2}\sum_{s=1}^{t-1}\sum_{k=0}^{s-1} u_t u_s
v_{t-k} v_{s-k}g_n(k)~\mbox{ and
}~W_{2nt}=\frac{2}{n^2}\sum_{s=2}^{t-1}\sum_{k=1}^{s-1}v_tv_s
u_{t-k}u_{s-k}g_n(-k).\] Then $W_{1nt}$ and $W_{2nt}$ form
martingale differences with respect to ${\cal F}_t$, where ${\cal
F}_t$ is the $\sigma$-field generated by
$\{u_s,v_s\}_{s=-\infty}^{t}$. Letting $W_{nt}=W_{1nt}+W_{2nt}$,
then $H_{11n}+H_{21n}=n^{-1}\sum_{t=3}^{n}W_{nt}+n^{-1}W_{1n2}$,
where the latter term is easily seen to be $o_p(1)$.

We shall apply the martingale central limit theorem  of Brown
(1971). It suffices to verify the following two conditions:
\begin{eqnarray}
\label{eq:conditionalvar}
\sigma_n^{-2}n^{-2}\sum_{t=3}^{n}\tilde{W}_{nt}^2\rightarrow_{p}
1,~~\mbox{where}~~\tilde{W}_{nt}^2=\E(W_{nt}^2|{\cal F}_{t-1})
\end{eqnarray}
\begin{eqnarray}
\label{eq:linderberg}
~\mbox{and}~~~\sigma_n^{-2}n^{-2}\sum_{t=3}^{n}\E[W_{nt}^2{\bf
1}\{|W_{nt}|>\epsilon n\sigma_n\}] \rightarrow 0,~\mbox{for
any}~\epsilon>0.
\end{eqnarray}
Since
$\sigma_n^{-2}n^{-2}\sum_{t=3}^{n}\E(\tilde{W}_{nt}^2)=\sigma_n^{-2}\var(H_{11n}+H_{21n})+o(1)=1+o(1)$,
as we have shown,  (\ref{eq:conditionalvar}) is implied by
\begin{eqnarray}
\label{eq:Wjnt}
\var\left(\sigma_n^{-2}n^{-2}\sum_{t=3}^n\tilde{W}_{jnt}^2\right)=o(1),~\mbox{where}~\tilde{W}_{jnt}^2=\E(W_{jnt}^2|{\cal
F}_{t-1}),~ j=1,2,
\end{eqnarray}
 which is established in Lemma~\ref{lem:varw1n}. To
prove (\ref{eq:linderberg}), it suffices in view of $(a+b)^4\le
8(a^4+b^4), a,b\in\R$ to verify
\[\sigma_n^{-4}n^{-4}\sum_{t=3}^{n}\E (W_{jnt}^4)=o(1),~j=1,2,\]
which can be shown by a similar and slightly simpler argument as
in the proof of Lemma~\ref{lem:varw1n}. The details are omitted.

  Thus (\ref{eq:1121}) is established and the conclusion follows.

\qed
\begin{lemma}
\label{lem:varw1n} Under the assumptions of
Theorem~\ref{th:mainresult1}, the random variable
$\tilde{W}_{jnt}^2$ defined in (\ref{eq:Wjnt}) satisfies
\[{\rm var}\left(\sigma_n^{-2}n^{-2}\sum_{t=3}^n\tilde{W}_{jnt}^2\right)=o(1), ~j=1,2.\]
\end{lemma}

\noindent Proof of Lemma~\ref{lem:varw1n}: We shall only show the
case $j=1$, as the treatment for $j=2$ is similar. Note that
\begin{eqnarray*}
\E (W_{1nt}^2|{\cal
F}_{t-1})&=&\frac{4}{n^{4}}\E\left[\left(\sum_{s=2}^{t-1}\sum_{k=1}^{s-1}u_tu_s
v_{t-k}v_{s-k}g_n(k)+u_tv_t\sum_{s=1}^{t-1}u_s
v_sg_n(0)\right)^2|{\cal
F}_{t-1}\right]\\
&=& \frac{4}{n^{4}}\left[\left(\sum_{s=2}^{t-1}\sum_{k=1}^{s-1}u_s
v_{t-k}v_{s-k}g_n(k)\right)^2+\left(\sum_{s=1}^{t-1}u_s
v_sg_n(0)\right)^2\right].
\end{eqnarray*}
 Let $J_0:=\var(\sum_{t=3}^{n}(\sum_{s=1}^{t-1}u_s
v_sg_n(0))^2)$ and $J_1:=\var(\sum_{t=3}^n
(\sum_{s=2}^{t-1}\sum_{k=1}^{s-1}u_s v_{t-k}$ $v_{s-k}g_n(k))^2)$.
The conclusion follows from $J_k=o(\sigma_n^4n^{12})$, $k=0,1$. We
 consider $J_1$ first.
\begin{eqnarray*}
J_1&=&\sum_{t_1,t_3=3}^{n}\cov\left(\sum_{s_1,s_2=2}^{t_1-1}\sum_{k_1=1}^{s_1-1}\sum_{k_2=1}^{s_2-1}u_{s_1}u_{s_2}v_{t_1-k_1}v_{s_1-k_1}v_{t_1-k_2}v_{s_2-k_2}g_n(k_1)g_n(k_2),\right.\\
&&\left.\sum_{s_3,s_4=2}^{t_3-1}\sum_{k_3=1}^{s_3-1}\sum_{k_4=1}^{s_4-1}u_{s_3}
u_{s_4}v_{t_3-k_3}v_{s_3-k_3}v_{t_3-k_4}v_{s_4-k_4}g_n(k_3)g_n(k_4) \right)\\
&=&\sum_{t_1,t_3=3}^{n}\sum_{s_1,s_2=2}^{t_1-1}\sum_{k_1=1}^{s_1-1}\sum_{k_2=1}^{s_2-1}\sum_{s_3,s_4=2}^{t_3-1}\sum_{k_3=1}^{s_3-1}\sum_{k_4=1}^{s_4-1}g_n(k_1)g_n(k_2)g_n(k_3)g_n(k_4)J_{11},~\mbox{where}
\end{eqnarray*}
\begin{eqnarray*}
J_{11}&=&\cov(u_{s_1}u_{s_2}v_{t_1-k_1}v_{s_1-k_1}v_{t_1-k_2}v_{s_2-k_2},u_{s_3}
u_{s_4}v_{t_3-k_3}v_{s_3-k_3}v_{t_3-k_4}v_{s_4-k_4})\\
&=&\E(u_{s_1}u_{s_2}u_{s_3}
u_{s_4})\E(v_{t_1-k_1}v_{s_1-k_1}v_{t_1-k_2}v_{s_2-k_2}v_{t_3-k_3}v_{s_3-k_3}v_{t_3-k_4}v_{s_4-k_4})\\
&&-\E(u_{s_1}u_{s_2})\E(v_{t_1-k_1}v_{s_1-k_1}v_{t_1-k_2}v_{s_2-k_2})\E(u_{s_3}
u_{s_4})\\
&&\times\E(v_{t_3-k_3}v_{s_3-k_3}v_{t_3-k_4}v_{s_4-k_4})=J_{111}\times
J_{112}-J_{113}.
\end{eqnarray*}
In the above expression,
 $J_{111}=\E(u_{s_1}u_{s_2}u_{s_3} u_{s_4})={\bf
1}(s_1=s_2,s_3=s_4)+{\bf 1}(s_1=s_3,s_2=s_4)+{\bf
1}(s_1=s_4,s_2=s_3)+(c_4(u)-3){\bf 1}(s_1=s_2=s_3=s_4)$,
$J_{113}={\bf 1}(s_1=s_2,s_3=s_4,k_1=k_2,k_3=k_4)$ and
\begin{eqnarray*}
J_{112}&=&\E(v_{t_1-k_1}v_{s_1-k_1}v_{t_1-k_2}v_{s_2-k_2}v_{t_3-k_3}v_{s_3-k_3}v_{t_3-k_4}v_{s_4-k_4})\\
&=&\sum_{g}\cum(v_{i_j},i_j\in g_1)\cdots\cum(v_{i_j},i_j\in g_p),
\end{eqnarray*}
where $\sum_g$ is over all  partitions $g=\{g_1\cup
\cdots \cup g_p\}$ of the index set $\{t_1-k_1,t_1-k_2,t_3-k_3,t_3-k_4,s_1-k_1,s_2-k_2,s_3-k_3,s_4-k_4\}$.
Since $\E(v_t)=0$, only partitions $g$ with $\#g_i>1$ for all $i$
contribute.
 We shall
divide all contributing partitions into the following several
types and treat them one by one.
\begin{enumerate}
\item $\#g_1=\#g_2=4$. A typical term that contributes is
\[\cum(v_{t_1-k_1},v_{t_1-k_2},
       v_{t_3-k_3},v_{t_3-k_4})\cum(v_{s_1-k_1},v_{s_2-k_2},v_{s_3-k_3},v_{s_4-k_4}),\]
which equals to ${\bf
1}(k_1=k_2,k_3=k_4,s_1=s_2,s_3=s_4,t_1-k_1=t_3-k_3,s_1-k_1=s_3-k_3)$.

 \item $\#g_1=\#g_2=3,\#g_3=2$. A typical term is
 \[\cum(v_{t_1-k_1},v_{t_1-k_2},v_{t_3-k_3})\cum(v_{s_1-k_1},v_{s_2-k_2},v_{t_3-k_4})\cov(v_{s_3-k_3},v_{s_4-k_4}),\]
which is ${\bf
1}(k_1=k_2,s_1=s_2,t_1-k_1=t_3-k_3,t_3-k_4=s_1-k_1,s_3-k_3=s_4-k_4)$.

 \item $\#g_1=\#g_2=\#g_3=\#g_4=2$. One such term is
\begin{eqnarray}
\label{eq:2222}\cov(v_{t_1-k_1},v_{t_1-k_2})\cov(v_{t_3-k_3},v_{t_3-k_4})\cov(v_{s_1-k_1},v_{s_2-k_2})\cov(v_{s_3-k_3},v_{s_4-k_4}),
\end{eqnarray}
which is ${\bf 1}(k_1=k_2,k_3=k_4,s_1=s_2,s_3=s_4)$. Note that
(\ref{eq:2222}) multiplied with ${\bf 1}(s_1=s_2,s_3=s_4)$ in
$J_{111}$ cancels out $J_{113}$. Thus all non-vanishing terms in
$J_{11}$ involve at least five restrictions on the indices
$t_1,t_3,s_1,s_2,s_3,s_4,k_1,k_2,k_3,k_4$.
\end{enumerate}
In the following we will find a bound for only one such term since
the derivation for other terms is similar. For example, one of the
terms in $J_{1}$, denoted by $E^*$, corresponds to the case
$k_1=k_2$, $k_3=k_4$, $s_1=s_2=s_3=s_4$, i.e.
\begin{eqnarray*}
E^*&=&\sum_{t_1,t_3=3}^{n}\sum_{s_1=2}^{(t_1-1)\wedge(t_3-1)}\sum_{k_1,k_3=1}^{s_1-1}g_n^2(k_1)
g_n^2(k_3).
\end{eqnarray*}
 Recall
$g_n(k)=\sum_{j,j'=1}^{n-1}G_{j,j'}\Lambda_n(\lambda_{j-j'})e^{ik\lambda_{j-j'}}$.
We have
\begin{eqnarray*}
&&\hspace{-0.5cm}E^*=\sum_{j_1,j_1',j_2,j_2'=1}^{n-1}G_{j_1,j_1'}G_{j_2,j_2'}\sum_{j_3,j_3',j_4,j_4'=1}^{n-1}G_{j_3,j_3'}G_{j_4,j_4'}\Lambda_n(\lambda_{j_1-j_1'})\Lambda_n(\lambda_{j_2-j_2'})\\
&&\hspace{0.05cm}\Lambda_n(\lambda_{j_3-j_3'})\Lambda_n(\lambda_{j_4-j_4'})\sum_{t_1,t_3=3}^{n}\sum_{s_1=2}^{(t_1-1)\wedge(t_3-1)}
\sum_{k_1,k_3=1}^{s_1-1}e^{ik_1\lambda_{j_1-j_1'+j_2-j_2'}}e^{ik_3\lambda_{j_3-j_3'+j_4-j_4'}}.
\end{eqnarray*}
Then
\begin{eqnarray}
\label{eq:boundE}
 &&|E^*|\le
Cn^4\sum_{m_1,m_2,m_3,m_4=2-n}^{n-2}|\Lambda_n(\lambda_{m_1})||\Lambda_n(\lambda_{m_2})|
|\Lambda_n(\lambda_{m_3})||\Lambda_n(\lambda_{m_4})|\\
&&\hspace{0.5cm}\times\left|\sum_{t_1,t_3=3}^{n}\sum_{s_1=2}^{(t_1-1)\wedge(t_3-1)}
\sum_{k_1,k_3=1}^{s_1-1}e^{ik_1\lambda_{m_1+m_2}}e^{ik_3\lambda_{m_3+m_4}}\right|
~(=: E^*_{m_1m_2m_3m_4}).\nonumber
\end{eqnarray}
It is not hard to see that $E^*_{m_1m_2m_3m_4}$
 is $O(|\lambda_{m_1+m_2}|^{-1}|\lambda_{m_3+m_4}|^{-1}n^3)$ if
$m_1+m_2\notin \{0,\pm n\}$ and $m_3+m_4\not\in \{0,\pm n\}$;
$O(|\lambda_{m_3+m_4}|^{-1}n^{4})$ if $m_1+m_2\in \{0,\pm n\}$ and
$m_3+m_4\notin \{0,\pm n\}$; $O(|\lambda_{m_1+m_2}|^{-1}n^{4})$ if
$m_1+m_2\notin \{0,\pm n\}$ and $m_3+m_4\in \{0,\pm n\}$; $O(n^5)$
if $m_1+m_2\in \{0,\pm n\}$ and $m_3+m_4\in \{0,\pm n\}$.


Combined with (\ref{eq:Lambdan}), we can derive from
(\ref{eq:boundE}) that $|E^*|=O(n^{11}B_n^2\log^4
n)=o(n^{12}B_n^2)$. Therefore $|J_{1}|=o(\sigma_n^4 n^{12})$ and a
similar argument yields $J_0=o(\sigma_n^4 n^{12})$. The conclusion
follows.

\qed

\begin{lemma}
\label{lem:variance}
 Under the assumptions of Theorem~\ref{th:mainresult1}, the random
 variable $H_n$ in (\ref{eq:Hn}) satisfies
\[\E(H_n)=B_n s(K)+O(1)~~\mbox{and}~~{\rm var}(H_n)=2 B_n d(K)+o(B_n).\]
\end{lemma}
\noindent Proof of Lemma~\ref{lem:variance}: Recall from
Lemma~\ref{lem:clt} that $H_n=4\pi^2
n^{-1}\sum_{j,j'=1}^{n-1}\Lambda_n(\lambda_{j-j'})$
$G_{j,j'}w_{uj}\overline{w_{vj}} \cdot \overline{w_{uj'}}
w_{vj'}$. Since
$\E(w_{uj}\overline{w_{uj'}})=\E(\overline{w_{vj}}w_{vj'})=(2\pi)^{-1}{\bf
1}(j-j'=0~\mbox{mod}~n)$, we have
$\E(H_n)=(1-1/n)\sum_{h_1=-B_n}^{B_n}K^2(h_1 b_n)=B_ns(K)+O(1)$
under Assumption~\ref{as:kernel}.

We proceed to calculate $\var(H_n)$. Note that
\begin{eqnarray*}
\mbox{var}(H_n)
&=&\frac{16\pi^4}{n^2}\sum_{j_1,j_1'=1}^{n-1}\sum_{j_2,j_2'=1}^{n-1}G_{j_1,j_1'}\overline{G_{j_2,j_2'}}\Lambda_n(\lambda_{j_1-j_1'})\Lambda_n(\lambda_{j_2-j_2'})C(j_1,j_2,j_1',j_2'),
\end{eqnarray*}
where by Lemma~\ref{lem:wwww},
\begin{eqnarray}
\label{eq:Cjjjj}
&&\hspace{-0.5cm}C(j_1,j_2,j_1',j_2'):=\mbox{cov}(w_{uj_1}\overline{w_{vj_1}}
\cdot \overline{w_{uj_1'}}w_{vj_1'},\overline{w_{uj_2}}{w_{vj_2}}
\cdot {w_{uj_2'}}\overline{w_{vj_2'}})\\
&&\hspace{2.1cm}=\E(w_{uj_1}\overline{w_{u j_1'}} \overline{w_{u
j_2}}{w_{u j_2'}}) \E(\overline{w_{v j_1}}w_{v
j_1'} {w_{vj_2}} \overline{w_{v j_2'}})\nonumber\\
&&\hspace{0.0cm}-(16\pi^4)^{-1}{\bf 1}(j_1=j_1'){\bf
1}(j_2=j_2')=(16\pi^4)^{-1}[{\bf 1}(j_1+j_2'=n){\bf
1}(j_1'+j_2=n)\nonumber\\
&&\hspace{0.0cm}+{\bf 1}(j_1=j_2){\bf
1}(j_2'=j_1')+(c_4(u)-3)(c_4(v)-3)n^{-2}{\bf
1}(j_1+j_2'-j_1'-j_2=0,\pm
n)\nonumber\\
&&\hspace{0.0cm}+2{\bf 1}(j_1+j_2'=n){\bf 1}(j_1=j_1'){\bf
1}(j_2=j_2')+2{\bf 1}(j_1+j_2'=n){\bf 1}(j_1=j_2){\bf
1}(j_2'=j_1')\nonumber\\
&&\hspace{0.0cm}+2{\bf
1}(j_1=j_2=j_1'=j_2')+(c_4(u)+c_4(v)-6)n^{-1}\{{\bf
1}(j_1=j_1'){\bf
1}(j_2=j_2')\nonumber\\
&&\hspace{0.0cm}+{\bf 1}(j_1+j_2'=n){\bf 1}(j_1'+j_2=n)+{\bf
1}(j_1=j_2){\bf 1}(j_2'=j_1')\} ].\nonumber
\end{eqnarray}
 A simple calculation shows that only the first two terms in
 $C(j_1,j_2,j_1',j_2')$, i.e. $(16\pi^4)^{-1}[{\bf 1}(j_1+j_2'=n){\bf
1}(j_1'+j_2=n)+{\bf 1}(j_1=j_2){\bf 1}(j_2'=j_1')]$ contribute to
the dominant term. Therefore,
\begin{eqnarray*}
\mbox{var}(H_n)&=&\frac{2(1+o(1))}{n^2}\sum_{j_1,j_1'=1}^{n-1}\Lambda_n^2(\lambda_{j_1-j_1'})\\
&=&\frac{2(1+o(1))}{n^2}\sum_{j_1,j_1'=1}^{n-1}\sum_{h_1,h_3=-B_n}^{B_n}K^2(h_1b_n)
K^2(h_3b_n) e^{i(h_1-h_3)\lambda_{j_1-j_1'}}\\
&=&\frac{4(1+o(1))}{n^2}\sum_{h_1,h_3=-B_n}^{B_n}K^2(h_1b_n)
K^2(h_3b_n)
\sum_{k=1}^{n}k\cos(k\lambda_{h_1-h_3})\\
&=&\frac{4(1+o(1))}{n^2}\sum_{k=1}^{n}k [a_n^2(k)
 +b_n^2(k)],
\end{eqnarray*}
where $a_n(k)=\sum_{h=-B_n}^{B_n}K^2(hb_n)\cos(k\lambda_h)$ and
$b_n(k)=\sum_{h=-B_n}^{B_n}K^2(hb_n)\sin(k\lambda_h)$.

By Lemma 4 in Wu and Shao (2007),
\[\sum_{k=1}^{n}ka_n^2(k)=\frac{n^2}{4}\left[\sum_{h=1}^{B_n}4K^4(hb_n)+1\right]+O(nB_n^2).\]
Since $K(\cdot)$ is symmetric, $b_n(k)=0$. Thus
\[\var(H_n)=4\sum_{h=1}^{B_n}K^4(hb_n)+o(B_n)=2B_n d(K)+o(B_n).\]
This completes the proof.

 \qed

\noindent Proof of Theorem~\ref{th:mainresult2}: The proof
follows the argument in the proof of Chen and Deo's (2004) Theorem
4.
Let
${\pmb{\theta}}_{0}=({\pmb{\theta}}_{10}',{\pmb{\theta}}_{20}')'$
and
$\hat{{\pmb{\theta}}}=(\hat{{\pmb{\theta}}}_1',\hat{{\pmb{\theta}}}_2')'$.
It suffices to show that
$n(T_n-T_n(\hat{{\pmb{\theta}}}))=o_p(\sqrt{B_n})$, which follows
from the following two assertions:
\begin{eqnarray}
\label{eq:negdiff1}
\sum_{l=0}^{n-1}(|\tilde{f}_{X_1X_2}(\lambda_l)|^2-|\hat{f}_{X_1X_2}(\lambda_l)|^2)&=&o_p(\sqrt{B_n}),\\
\label{eq:negdiff2}
\frac{2\pi}{n}\sum_{l=0}^{n-1}\{\tilde{f}_{X_kX_k}(\lambda_l)-\hat{f}_{X_kX_k}(\lambda_l)\}&=&o_p(B_n^{-1/2}),
~k=1,2.
\end{eqnarray}
We shall only provide a proof for (\ref{eq:negdiff1}), since the
treatment of (\ref{eq:negdiff2}) is similar. For
${\pmb{\theta}}=({\pmb{\theta}}_1',{\pmb{\theta}}_2')'$ let
$G(\lambda_j,\lambda_{j'};{{{\pmb{\theta}}}})=\Pi_{k=1}^2f_k^{-1/2}(\lambda_j;{{\pmb{\theta}}}_k)f_k^{-1/2}(\lambda_{j'};{{\pmb{\theta}}}_k)$.
  Further let
$I_{12j}^*=I_{12j}/\sqrt{f_1(\lambda_j;{\pmb{\theta}}_{10})
f_2(\lambda_j;{\pmb{\theta}}_{20})}$. Then the LHS of
 (\ref{eq:negdiff1}) is
\begin{eqnarray*}
 &&\hspace{-0.5cm}\frac{4\pi^2}{n^2}\sum_{l=0}^{n-1}\sum_{j,j'=1}^{n-1}W(\lambda_{l-j})W(\lambda_{l-j'})I_{12j}\overline{I_{12j'}}\{G(\lambda_j,\lambda_{j'};{\hat{{\pmb{\theta}}}})-G(\lambda_j,\lambda_{j'};{\pmb{\theta}}_0)\}\\
 &&\hspace{0.5cm}=\frac{1}{n}\sum_{j,j'=1}^{n-1}\Lambda_n(\lambda_{j-j'})I_{12j}^*\overline{I_{12j'}^*}\frac{\{G(\lambda_j,\lambda_{j'};{\hat{{\pmb{\theta}}}})-G(\lambda_j,\lambda_{j'};{{\pmb{\theta}}_0})\}}{G(\lambda_j,\lambda_{j'};\pmb{\theta}_0)},
\end{eqnarray*}
where we have applied (\ref{eq:WLambda}). Let $I_1$ and $I_2$ be
the index set of $\Theta_1\subset\R^{q_1}$ and $\Theta_2\subset
\R^{q_2}$ respectively. For every $\lambda_j$ and $\lambda_{j'}$,
a Taylor series expansion yields
\begin{eqnarray*}
&&G(\lambda_j,\lambda_{j'};{\hat{{\pmb{\theta}}}})-G(\lambda_j,\lambda_{j'};{{\pmb{\theta}}_0})=\sum_{u\in
I_1}(\hat{\theta}_{1_u}-{\theta}_{10_u})\frac{\partial
G(\lambda_j,\lambda_{j'};{{\pmb{\theta}}_0})}{\partial
\theta_{1_u}}\\
&&\hspace{0.5cm}+\sum_{u\in
I_2}(\hat{\theta}_{2_u}-{\theta}_{20_u})\frac{\partial
G(\lambda_j,\lambda_{j'};{{\pmb{\theta}}_0})}{\partial
\theta_{2_u}}+\frac{1}{2}(\hat{{\pmb{\theta}}}-{\pmb{\theta}}_0)'\frac{\partial^2
G(\lambda_j,\lambda_{j'};\tilde{{\pmb{\theta}}}_{j j'})}{\partial
{\pmb{\theta}}^2}(\hat{{\pmb{\theta}}}-{\pmb{\theta}}_0),
\end{eqnarray*}
where
$\tilde{{\pmb{\theta}}}_{jj'}=(\tilde{{\pmb{\theta}}}_{1jj'}',\tilde{{\pmb{\theta}}}_{2jj'}')'={\pmb{\theta}}_0+\alpha_{j
j'}(\hat{{\pmb{\theta}}}-{\pmb{\theta}}_0)$ for some
$\alpha_{jj'}\in [0,1]$.

For $k=1,2$, let $g_{ku}(\lambda;{\pmb{\theta}}_k)=\partial \log
f_k(\lambda;{\pmb{\theta}}_k)/\partial{{\theta}}_{k_u}$ and
$g_{kuv}(\lambda;{\pmb{\theta}}_k)=\partial^2
/\partial{{\theta}}_{k_u}\partial\theta_{k_v}$ $\log
f_k(\lambda;{\pmb{\theta}}_{k})$,
 $u,v\in I_k$. Then for $u\in I_k$,
\begin{eqnarray}
\label{eq:1stder}
 \frac{\partial
G(\lambda_j,\lambda_{j'};{{\pmb{\theta}}_0})}{\partial
\theta_{k_u}}=-\frac{1}{2}
G(\lambda_j,\lambda_{j'};{\pmb{\theta}}_{0})[g_{ku}(\lambda_j;{\pmb{\theta}}_{k0})+g_{ku}(\lambda_{j'};{\pmb{\theta}}_{k0})].
\end{eqnarray}
Let $A_{u v}(\lambda_j,\lambda_{j'};\tilde{{\pmb{\theta}}}_{jj'})$
be the $(u, v)$th element of the matrix $\partial^2
G(\lambda_j,\lambda_{j'};{\tilde{{\pmb{\theta}}}_{jj'}})/\partial
{\pmb{\theta}}^2$. Then
\begin{eqnarray*}
&&\hspace{-1cm}A_{u
v}(\lambda_j,\lambda_{j'};\tilde{{\pmb{\theta}}}_{jj'})=\frac14
G(\lambda_j,\lambda_{j'};\tilde{{\pmb{\theta}}}_{jj'})
\{[g_{ku}(\lambda_j;\tilde{{\pmb{\theta}}}_{kjj'})+
g_{ku}(\lambda_{j'};\tilde{{\pmb{\theta}}}_{kjj'})]^2\\
&&-2[g_{kuv}(\lambda_j;\tilde{{\pmb{\theta}}}_{kjj'})+g_{kuv}(\lambda_{j'};\tilde{{\pmb{\theta}}}_{kjj'})]\}
,~u,v\in I_k, k=1,2.\\
&&=\frac14
G(\lambda_j,\lambda_{j'};\tilde{{\pmb{\theta}}}_{jj'})\{
[g_{k_1u}(\lambda_j;\tilde{{\pmb{\theta}}}_{k_1
jj'})+g_{k_1u}(\lambda_{j'};\tilde{{\pmb{\theta}}}_{k_1
jj'})][g_{k_2v}(\lambda_j; \tilde{{\pmb{\theta}}}_{k_2
jj'})\\
&&+g_{k_2v}(\lambda_{j'};\tilde{{\pmb{\theta}}}_{k_2 jj'})]\},
~u\in I_{k_1}, v\in I_{k_2}, (k_1,k_2)=(1,2)~\mbox{or}~(2,1).
\end{eqnarray*}
Since
$|\hat{{\pmb{\theta}}}_k-{\pmb{\theta}}_{k0}|=O_{p}(n^{-1/2})$ and
(\ref{eq:1stder}), (\ref{eq:negdiff1}) follows from the following
two statements:
\begin{enumerate}
\item For each $u\in I_k$, $k=1,2$,
\begin{eqnarray}
\label{eq:statement1}
\hspace{-0.5cm}L_n:=\sum_{j,j'=1}^{n-1}\Lambda_n(\lambda_{j-j'})I_{12j}^*\overline{I_{12j'}^*}[g_{ku}(\lambda_j;{\pmb{\theta}}_{k0})+g_{ku}(\lambda_{j'};{\pmb{\theta}}_{k0})]=o_p(n^{3/2}\sqrt{B_n}).
\end{eqnarray}
\item For every $(u,v)$,
\begin{eqnarray}
\label{eq:statement2}
\hspace{-0.5cm}M_n:=\sum_{j,j'=1}^{n-1}\Lambda_n(\lambda_{j-j'})I_{12j}^*\overline{I_{12j'}^*}\frac{A_{u
v}(\lambda_j,\lambda_{j'};\tilde{{\pmb{\theta}}}_{jj'})}{G(\lambda_j,\lambda_{j'};\pmb{\theta}_0)}=o_p(n^2\sqrt{B_n}).
\end{eqnarray}
\end{enumerate}
By (\ref{eq:bartlettX}) and (\ref{eq:bartlettY}), $L_{n}$
 can be written as a sum of 16 terms, of which the dominant term
is given by
\begin{eqnarray*}
L_{1n}:=4\pi^2\sum_{j,j'=1}^{n-1}\Lambda_n(\lambda_{j-j'})G_{j,j'}w_{uj}\overline{w_{vj}}\cdot\overline{w_{uj'}}w_{vj'}[g_{ku}(\lambda_j;{\pmb{\theta}}_{k0})+g_{ku}(\lambda_{j'};{\pmb{\theta}}_{k0})].
\end{eqnarray*}
Then $\E|L_{1n}|^2$ is
\begin{eqnarray*}
&&\hspace{-0.5cm}16\pi^4\sum_{j_1,j_1'=1}^{n-1}\sum_{j_2,j_2'=1}^{n-1}\Lambda_n(\lambda_{j_1-j_1'})\Lambda_n(\lambda_{j_2-j_2'}) G_{j_1,j_1'}\overline{G_{j_2,j_2'}}[g_{ku}(\lambda_{j_1};{\pmb{\theta}}_{k0})+g_{ku}(\lambda_{j_1'};{\pmb{\theta}}_{k0})]\\
&&\hspace{0.5cm}\times[g_{ku}(\lambda_{j_2};{\pmb{\theta}}_{k0})+g_{ku}(\lambda_{j_2'};{\pmb{\theta}}_{k0})]\E(w_{uj_1}\overline{w_{uj_1'}}\cdot
\overline{w_{uj_2}}w_{uj_2'})\E(\overline{w_{vj_1}}
w_{vj_1'}w_{vj_2}\overline{w_{vj_2'}}).
\end{eqnarray*}
In view of Lemma~\ref{lem:wwww} and (\ref{eq:Cjjjj}), the dominant
term of  $\E |L_{1n}|^2$ is a sum of three terms, which correspond
to $j_1=j_1'$ and $j_2=j_2'$, $j_1+j_2'=n$ and $j_1'+j_2=n$ as
well as  $j_1=j_2$ and $j_2'=j_1'$. Since
$|g_{ku}(\lambda;{\pmb{\theta}}_{k0})|\le C|\lambda|^{-\delta}$
 under Assumption~\ref{as:condiff}, all the three terms can be
bounded by
\[CB_n^2\sum_{j_1,j_2=1}^{n-1}(\lambda_{j_1}\lambda_{j_2})^{-2\delta}=O(n^2 B_n^2), ~~\mbox{for}~\delta\in (0,1/2).\]


 For the other terms involved in $\E|L_{1n}|^2$,  a tighter bound
 than $O(n^2B_n^2)$ can be established using Lemma~\ref{lem:bound1} and
(\ref{eq:Lambdan}). So $L_{1n}=O_p(nB_n)=o_{p}(n^{3/2}B_n^{1/2})$.
The remaining 15 terms in $L_{n}$ can be shown to be
$o_p(n^{3/2}B_n^{1/2})$ using the bounds established in
Lemma~\ref{lem:bound1}.  Since the proof does not involve
additional methodological difficulties, we omit the details.

 Now we prove
(\ref{eq:statement2}). Note that
 $\tilde{{\pmb{\theta}}}_{kjj'}=(\tilde{{\pmb{\beta}}}_{kjj'}',\tilde{d}_{kjj'})'$,
 $k=1,2$. Under Assumptions~\ref{as:longmemo} and~\ref{as:condiff}, we have
\begin{eqnarray*}
\label{eq:boundA} |A_{u
v}(\lambda_j,\lambda_{j'};\tilde{{\pmb{\theta}}}_{jj'})|\le C
(\lambda_j\lambda_{j'})^{\tilde{d}_{1jj'}+\tilde{d}_{2jj'}}(\lambda_j^{-\delta}+\lambda_{j'}^{-\delta})^2
\end{eqnarray*}
uniformly in $(u,v)$. By  Lemma~\ref{lem:bound1} and
(\ref{eq:bartlettX}), uniformly in $j=1,2,\cdots,n-1$,
\begin{eqnarray*}
\label{eq:boundX}
 \frac{\E
\{I_{X_1X_1}(\lambda_j)\}}{f_{1}(\lambda_j;\pmb{\theta}_{10})}\le
2 [2\pi \E|w_{uj}|^2+\E|R_{uj}|^2]\le C.
\end{eqnarray*}
 Hence, by the Cauchy-Schwarz inequality, uniformly in $j,j'=1,2,\cdots,n-1$,
\begin{eqnarray}
\label{eq:boundperiod}\E|I_{12j}^*\overline{I_{12j'}^*}|\le
\sup_{j=1,\cdots,n-1}\frac{\E\{
I_{X_1X_1}(\lambda_j)\}}{f_1(\lambda_j;{\pmb{\theta}}_{10})}
\sup_{j=1,\cdots,n-1}\frac{\E\{
I_{X_2X_2}(\lambda_j)\}}{f_2(\lambda_j;{\pmb{\theta}}_{20})}<C.
\end{eqnarray}
In the sequel, we shall treat the following several cases
separately.
\begin{enumerate}
\item $\hat{d}_1\ge d_{10}$ and $\hat{d}_{2}\ge d_{20}$. Then
$\tilde{d}_{1jj'}\ge d_{10}$ and $\tilde{d}_{2jj'}\ge d_{20}$ for
all $j,j'$, which implies that
$|G^{-1}(\lambda_j,\lambda_{j'};{\pmb{\theta}}_0)A_{u
v}(\lambda_j,\lambda_{j'};\tilde{{\pmb{\theta}}}_{jj'})|\le
C(\lambda_{j}^{-\delta}+\lambda_{j'}^{-\delta})^2$. In this case,
we get
\begin{eqnarray*}
|M_n|{\bf 1}(\hat{d}_1\ge d_{10},\hat{d}_{2}\ge d_{20})
     \le C\sum_{j,j'=1}^{n-1}|I_{12j}^*\overline{I_{12j'}^*}||\Lambda_n(\lambda_{j-j'})|(\lambda_{j}^{-\delta}+\lambda_{j'}^{-\delta})^2.
\end{eqnarray*}
Then by (\ref{eq:boundperiod}) and (\ref{eq:Lambdan}),
$\E|M_n|{\bf 1}(\hat{d}_1\ge d_{10},\hat{d}_{2}\ge d_{20})$ is
bounded by
\[ C\sum_{j,j'=1}^{n-1}|\Lambda_n(\lambda_{j-j'})| (\lambda_{j}^{-\delta}+\lambda_{j'}^{-\delta})^{2}=O(n^2\log
 n).\]
  So $|M_n|{\bf 1}(\hat{d}_1\ge d_{10},\hat{d}_{2}\ge
d_{20})=o_p(n^2\sqrt{B_n})$ under Assumption~\ref{as:bandwidth}.

\item  $\hat{d}_1< d_{10}$ and $\hat{d}_2< d_{20}$. Denote by
$\Delta_{k}^{jj'}(\lambda_j)=f_k^{1/2}(\lambda_j;{\pmb{\theta}}_{k0})
f_k^{-1/2}(\lambda_j;\tilde{{\pmb{\theta}}}_{kjj'})-1$, $k=1,2$.
Under Assumption~\ref{as:longmemo} and (\ref{eq:lipschitz}), for
$j=1,2,\cdots,n-1$,
\begin{eqnarray}
\label{eq:Delta} |\Delta_{k}^{jj'}(\lambda_j)|\le
C|\hat{{\pmb{\theta}}}_k-{\pmb{\theta}}_{k0}|\lambda_j^{-(d_{k0}+\kappa)},~k=1,2.
\end{eqnarray}
 Let
$N_{j,j'}:=G^{-1}(\lambda_j,\lambda_{j'};\tilde{\pmb{\theta}}_{jj'})A_{u
v}(\lambda_j,\lambda_{j'};\tilde{{\pmb{\theta}}}_{jj'})$. Then
\begin{eqnarray}
\label{eq:boundMn}
 &&\hspace{-1cm}|M_n|{\bf 1}(\hat{d}_1< d_{10},\hat{d}_2<
d_{20})\le\sum_{j,j'=1}^{n-1}|\Lambda_n(\lambda_{j-j'})||I_{12j}^*\overline{I_{12j'}^*}||N_{j,j'}|\nonumber\\
&&\hspace{-0.5cm}\times(1+|\Delta_{1}^{jj'}(\lambda_j)|)(1+|\Delta_{1}^{jj'}(\lambda_{j'})|)(1+|\Delta_{2}^{jj'}(\lambda_j)|)(1+|\Delta_{2}^{jj'}(\lambda_{j'})|).
\end{eqnarray}
The RHS (right hand side) of (\ref{eq:boundMn}) consists of 16 terms, one of which
is
\[J_{n}:=\sum_{j,j'=1}^{n-1}|\Lambda_n(\lambda_{j-j'})||I_{12j}^*\overline{I_{12j'}^*}||\Delta_{1}^{jj'}(\lambda_j)||N_{j,j'}|.\]
 According to (\ref{eq:Delta}) and the fact that $|N_{j,j'}|\le
C(\lambda_j^{-\delta}+\lambda_{j'}^{-\delta})^2$, we have
\begin{eqnarray*}
|J_{n}|\le C|\hat{{\pmb{\theta}}}_1-{{\pmb{\theta}}}_{10}|
\left\{\sum_{j,j'=1}^{n-1}|\Lambda_n(\lambda_{j-j'})||I_{12j}^*\overline{I_{12j'}^*}|\lambda_j^{-(d_{10}+\kappa)}(\lambda_{j}^{-\delta}+\lambda_{j'}^{-\delta})^2\right\},
\end{eqnarray*}
where the expectation of the  bracketed term above is bounded by
\[C\sum_{j,j'=1}^{n-1}|\Lambda_n(\lambda_{j-j'})|\lambda_j^{-(d_{10}+\kappa)}(\lambda_{j}^{-\delta}+\lambda_{j'}^{-\delta})^2 =O(n^{2}\log n)\]
for $\delta\in (0,(1-d_{10}-\kappa) /2)$. Together with the fact
that $|\hat{{\pmb{\theta}}}_1-{\pmb{\theta}}_{10}|=O_p(n^{-1/2})$,
we get $J_n=o_p(n^2\sqrt{B_n})$. The other terms on the RHS of
(\ref{eq:boundMn}) can be shown to be $o_{p}(n^2\sqrt{B_n})$  in a
similar fashion. Hence $|M_n|{\bf 1}(\hat{d}_1< d_{10},\hat{d}_2<
d_{20})=o_p(n^2\sqrt{B_n})$.

 \item $\hat{d}_1\ge d_{10}, \hat{d}_2<
d_{20}$ and $\hat{d}_2\ge d_{20}, \hat{d}_1< d_{10}$. These two
cases can be handled in a similar manner as in the previous two
cases and the details are omitted.
\bigskip

The proof is now complete.
\end{enumerate}

\qed

\noindent Proof of Theorem~\ref{th:mainresult3}: Similar to the
proof of Theorem~\ref{th:mainresult2}, it suffices to show the
following two assertions:
\begin{eqnarray}
\label{eq:pneg1}
\sum_{l=0}^{n-1}(|\check{f}_{X_1X_2}(\lambda_l)|^2-|\hat{f}_{X_1X_2}(\lambda_l)|^2)&=&o_p(\sqrt{B_n}),\\
\label{eq:pneg2}
\frac{2\pi}{n}\sum_{l=0}^{n-1}\{\check{f}_{X_kX_k}(\lambda_l)-\hat{f}_{X_kX_k}(\lambda_l)\}&=&o_p(B_n^{-1/2}),~k=1,2.
\end{eqnarray}
We only prove (\ref{eq:pneg1}), as the proof of (\ref{eq:pneg2})
is similar. Let
\begin{eqnarray*}
H_{jj'}=\Pi_{k=1}^{2}f_{X_kX_k}^{-1/2}(\lambda_j)f_{X_kX_k}^{-1/2}(\lambda_{j'}),~
H_{jj'}({\pmb{\gamma}})=\Pi_{k=1}^{2}f_{k}^{-1/2}(\lambda_j;{\pmb{\gamma}}_k)f_{k}^{-1/2}(\lambda_{j'};{\pmb{\gamma}}_k).
\end{eqnarray*}
Denote by
${\pmb{\gamma}}_{0}=({\pmb{\gamma}}_{10}',{\pmb{\gamma}}_{20}')'$.
Applying (\ref{eq:WLambda}), the LHS of (\ref{eq:pneg1}) is
\begin{eqnarray*}
&&\frac{4\pi^2}{n^2}\sum_{l=0}^{n-1}\sum_{j,j'=1}^{n-1}W(\lambda_{l-j})W(\lambda_{l-j'})I_{12j}\overline{I_{12j'}}\{H_{jj'}(\hat{{\pmb{\gamma}}})-H_{jj'}\}\\
&&\hspace{-0.5cm}=\frac{1}{n}\sum_{j,j'=1}^{n-1}\Lambda_n(\lambda_{j-j'})I_{12j}\overline{I_{12j'}}\{[H_{jj'}({\pmb{\gamma}}_0)-H_{jj'}]+[H_{jj'}(\hat{{\pmb{\gamma}}})-H_{jj'}({\pmb{\gamma}}_{0})]\}=S_{1n}+S_{2n}.
\end{eqnarray*}
Note that Assumption~\ref{as:expondecay} implies that
$|H_{jj'}({\pmb{\gamma}}_0) H_{jj'}^{-1}-1|\le C
\sum_{k=1}^{2}\sum_{l=p_k+1}^{\infty}|a_{kl}^0|$.
 Therefore, by (\ref{eq:Lambdan}),
(\ref{eq:boundperiod}) and (\ref{eq:pB1}),
\begin{eqnarray*}
\E|S_{1n}|&\le&
n^{-1}\sum_{j,j'=1}^{n-1}|\Lambda_n(\lambda_{j-j'})|\E|I_{12j}^*
\overline{I_{12j'}^*}|
|H_{jj'}({\pmb{\gamma}}_0)H_{jj'}^{-1}-1|\\
&\le&Cn\log n\sum_{k=1}^2\sum_{l=p_k+1}^{\infty}|a_{kl}^0|
=o(\sqrt{B_n}).
\end{eqnarray*}
 Regarding $S_{2n}$, we have
 \[S_{2n}=n^{-1}\sum_{j,j'=1}^{n-1}\Lambda_n(\lambda_{j-j'})I^*_{12j}\overline{I_{12j'}^*}\frac{H_{jj'}({\pmb{\gamma}}_0)}{H_{jj'}}\frac{[H_{jj'}(\hat{{\pmb{\gamma}}})-H_{jj'}({\pmb{\gamma}}_0)]}{H_{jj'}({\pmb{\gamma}}_0)}.\]
The remaining proof largely follows the argument in the proof of
Theorem~\ref{th:mainresult2}. Here we only outline key steps. Let
$J_k=\{1,\cdots,p_k+1\}$ be the index set for ${\pmb{\gamma}}_k$.
For each $(j,j')$, we apply a Taylor's expansion and obtain
\begin{eqnarray}
\label{eq:taylor}
H_{jj'}(\hat{{\pmb{\gamma}}})-H_{jj'}({\pmb{\gamma}}_0)&=&\sum_{k=1}^{2}\sum_{u\in
J_k}(\hat{{{\gamma}}}_{k_u}-{\gamma}_{{k0}_u})\frac{\partial
H_{jj'}({\pmb{\gamma}}_0)}{\partial \gamma_{k_u}}\\
&&+\frac12(\hat{{\pmb{\gamma}}}-{\pmb{\gamma}}_0)'\frac{\partial^2
H_{jj'}(\tilde{{\pmb{\gamma}}}_{jj'})}{\partial{\pmb{\gamma}}^2}(\hat{{\pmb{\gamma}}}-{\pmb{\gamma}}_0),\nonumber
\end{eqnarray}
where
$\tilde{{\pmb{\gamma}}}_{jj'}=(\tilde{{\pmb{\gamma}}}_{1jj'}',\tilde{{\pmb{\gamma}}}_{2jj'}')'={\pmb{\gamma}}_0+\alpha_{jj'}(\hat{{\pmb{\gamma}}}-{\pmb{\gamma}}_0)$
for some $\alpha_{jj'}\in [0,1]$.

Let $h_{ku}(\lambda;{\pmb{\gamma}}_k)=\partial \log
f_{k}(\lambda;{\pmb{\gamma}}_k)/\partial\gamma_{k_u}$,
$h_{kuv}(\lambda;{\pmb{\gamma}}_k)=\partial^2 \log
f_{k}(\lambda;{\pmb{\gamma}}_k)/\partial\gamma_{k_u}\partial\gamma_{k_v}$
for $u,v\in J_k$ and
$B_{uv}(\lambda_j,\lambda_{j'};\tilde{{\pmb{\gamma}}}_{jj'})$ be
the $(u,v)$th element of the matrix
$\partial^2H_{jj'}(\tilde{{\pmb{\gamma}}}_{jj'})/\partial{\pmb{\gamma}}^2$.
Then we have
\[\frac{\partial H_{jj'}({{\pmb{\gamma}}}_0)}{\partial \gamma_{k_u}}=-\frac12 H_{jj'}({\pmb{\gamma}}_0)[h_{ku}(\lambda_j;{\pmb{\gamma}}_{k0})+h_{ku}(\lambda_{j'};{\pmb{\gamma}}_{k0})]\]
and
\begin{eqnarray*}
&&B_{uv}(\lambda_j,\lambda_{j'};\tilde{{\pmb{\gamma}}}_{jj'})=\frac14
H_{jj'}(\tilde{{\pmb{\gamma}}}_{jj'})\{[h_{ku}(\lambda_j;\tilde{{\pmb{\gamma}}}_{kjj'})+h_{ku}(\lambda_{j'};\tilde{{\pmb{\gamma}}}_{kjj'})]^2\\
&&\hspace{0.5cm}-2[h_{kuv}(\lambda_j;\tilde{{\pmb{\gamma}}}_{kjj'})+h_{kuv}(\lambda_{j'};\tilde{{\pmb{\gamma}}}_{kjj'})]\},
~u,v\in J_k, k=1,2.\\
&&\hspace{0.5cm}=\frac14
H_{jj'}(\tilde{{\pmb{\gamma}}}_{jj'})\{[h_{k_1u}(\lambda_j;\tilde{{\pmb{\gamma}}}_{k_1jj'})+h_{k_1u}(\lambda_{j'};\tilde{{\pmb{\gamma}}}_{k_1jj'})][h_{k_2
v}(\lambda_j;\tilde{{\pmb{\gamma}}}_{k_2jj'})\\
&&\hspace{0.5cm}+h_{k_2
v}(\lambda_{j'};\tilde{{\pmb{\gamma}}}_{k_2jj'})]\}, u\in J_{k_1},
v\in J_{k_2}, (k_1,k_2)=(1,2)~\mbox{or}~(2,1).
\end{eqnarray*}
Note that $H_{jj'}({\pmb{\gamma}}_0)H_{jj'}^{-1}\le C$ under
Assumption~\ref{as:expondecay}. To prove $S_{2n}=o_p(\sqrt{B_n})$,
it suffices in view of Assumption~\ref{as:consistent} and
(\ref{eq:taylor}) to show that
\begin{enumerate}
\item Uniformly in $u\in J_k$, $k=1,2$,
\begin{eqnarray*}
\hspace{-1cm}L_n^*:=\sum_{j,j'=1}^{n-1}\Lambda_n(\lambda_{j-j'})I_{12j}^*\overline{I_{12j'}^*}\frac{H_{jj'}({\pmb{\gamma}}_0)}{H_{jj'}}[h_{ku}(\lambda_j;{\pmb{\gamma}}_{k0})+h_{ku}(\lambda_{j'};{\pmb{\gamma}}_{k0})]=o_p(n^{3/2}\sqrt{B_n}/p_k).
\end{eqnarray*}
\item Uniformly in $u\in J_{k_1}$ and $v\in J_{k_2}$,
$k_1,k_2=1,2$,
\begin{eqnarray*}
M_n^*:=\sum_{j,j'=1}^{n-1}\Lambda_n(\lambda_{j-j'})I_{12j}^*\overline{I_{12j'}^*}\frac{H_{jj'}({\pmb{\gamma}}_0)}{H_{jj'}}\frac{B_{uv}(\lambda_j,\lambda_{j'};\tilde{{\pmb{\gamma}}}_{jj'})}{H_{jj'}({\pmb{\gamma}}_0)}=o_{p}(n^2\sqrt{B_n}/(p_{k_1}p_{k_2})).
\end{eqnarray*}
\end{enumerate}
It is easy to see that Assumption~\ref{as:expondecay} implies that
uniformly in $\lambda$, $(u,v)$ and ${\pmb{\gamma}}_k$,
\[|h_{ku}(\lambda;{\pmb{\gamma}}_k)|\le C|\lambda|^{-\delta},~|h_{kuv}(\lambda;{\pmb{\gamma}}_k)|\le C|\lambda|^{-\delta} ~\mbox{for any}~\delta>0.\]
 Similar to the treatment of $L_n$ in the proof of
Theorem~\ref{th:mainresult2}, we can derive $\E|L_n^*|^2=O(n^2
B_n^2)=o(n^3B_n/p_k^2)$ under (\ref{eq:pB2}). A probabilistic
bound for $M_n^*$ can be established in the same way as that for
$M_n$ (see (\ref{eq:statement2})).  Here we omit the details but
mention the following fact (\ref{eq:lipschitz2}), which is needed
in the proof. Under Assumption~\ref{as:expondecay}, we have
\begin{eqnarray}
\label{eq:lipschitz2}
|f_k^{1/2}(\lambda;{\pmb{\gamma}}_{k}^{(1)})-f_k^{1/2}(\lambda;{\pmb{\gamma}}_{k}^{(2)})|\le
C|{\pmb{\gamma}}_{k}^{(1)}-{\pmb{\gamma}}_k^{(2)}|f_{k}^{1/2}(\lambda;{\pmb{\gamma}}_{k}^{(2)})
\end{eqnarray}
uniformly for all $\lambda$ and all ${\pmb{\gamma}}_k^{(1)}$ and
${\pmb{\gamma}}_k^{(2)}$ such that $d_{k}^{(1)}<d_{k}^{(2)}$.

The conclusion is established.
 \qed

\begin{lemma}
\label{lem:wwww} For any $t_1,t_2,t_3,t_4\in\Z$,
$\E[u_{t_1}u_{t_2}u_{t_3}u_{t_4}]=(c_4(u)-3){\bf
1}(t_1=t_2=t_3=t_4)+{\bf 1}(t_1=t_2){\bf 1}(t_3=t_4)+{\bf
1}(t_1=t_3){\bf 1}(t_2=t_4)+{\bf 1}(t_1=t_4){\bf 1}(t_2=t_3)$.
Further, for $1\le j_1, j_1', j_2, j_2'\le n-1$,
$\E(w_{uj_1}\overline{w_{u j_1'}} w_{u j_2}\overline{w_{u
j_2'}})=(4\pi^2)^{-1}[{\bf 1}(j_1=j_1'){\bf 1}(j_2=j_2')+{\bf
1}(j_1+j_2=n){\bf 1}(j_1'+j_2'=n)+{\bf 1}(j_1=j_2'){\bf
1}(j_2=j_1')+(c_4(u)-3) n^{-1}{\bf 1}(j_1+j_2-j_1'-j_2'=0,\pm
n)]$.
\end{lemma}
\noindent Proof of Lemma~\ref{lem:wwww}: Since $\{u_t\}$ are mean
zero iid random variables, the first assertion follows from the
following fact:
\begin{eqnarray*}
\E[u_{t_1}u_{t_2}u_{t_3}u_{t_4}]&=&\cum(u_{t_1},u_{t_2},u_{t_3},u_{t_4})+\cov(u_{t_1},u_{t_2})\cov(u_{t_3},u_{t_4})\\
&&+\cov(u_{t_1},u_{t_3})\cov(u_{t_2},u_{t_4})
+\cov(u_{t_1},u_{t_4})\cov(u_{t_2},u_{t_3}).
\end{eqnarray*}
Applying the first assertion, we have
\begin{eqnarray*}
\label{eq:wwww} && \hspace{-0.5cm}\E(w_{uj_1}\overline{w_{u j_1'}}
w_{u j_2}\overline{w_{u j_2'}})=\frac{1}{4\pi^2
n^2}\sum_{t_1,t_2,t_3,t_4=1}^{n} \E[u_{t_1} u_{t_2} u_{t_3}
u_{t_4}] e^{i(t_1\lambda_{j_1}-t_2\lambda_{j_1'}+t_3\lambda_{j_2}-t_4\lambda_{j_2'})}\\
&&\hspace{0.5cm}=\frac{1}{4\pi^2
n^2}\left[\sum_{t_1,t_3=1}^{n}e^{i
t_1\lambda_{j_1-j_1'}}e^{it_3\lambda_{j_2-j_2'}}+\sum_{t_1,t_2=1}^{n}
e^{it_1\lambda_{j_1+j_2}} e^{-it_2\lambda_{j_1'+j_2'}}\right.\\
&&\hspace{0.5cm}\left.+\sum_{t_1,t_2=1}^{n}
e^{it_1\lambda_{j_1-j_2'}}
e^{it_2\lambda_{j_2-j_1'}}+(c_4(u)-3)\sum_{t_1=1}^{n}
e^{it_1\lambda_{j_1+j_2-j_1'-j_2'}}\right].\\
\end{eqnarray*}
Since $\sum_{t=1}^{n}e^{it\lambda_j}=n{\bf 1}(j=0 ~\mbox{mod}~
n)$, the second assertion follows.

\qed

We now provide an auxiliary lemma on the bound of the second and
fourth cumulants of $R_{uj}$.
\begin{lemma}
\label{lem:bound1} Under Assumption~\ref{as:dependence}, we have
that
\begin{eqnarray*}
\label{eq:rrrr} |{\rm cov}(w_{uj_1},R_{uj_2})|&=&O(|j_2|^{-1/2})\\
\label{eq:rrrr1} |{\rm cov}(R_{uj_1},R_{uj_2})|&=&O(|j_1|^{-1/2}|j_2|^{-1/2})\\
\label{eq:rrrr2}
|\E(R_{uj_1}R_{uj_2}R_{uj_3}R_{uj_4})|&=&O(|j_1|^{-1/2}|j_2|^{-1/2}|j_3|^{-1/2}|j_4|^{-1/2})\\
\label{eq:rrrr3} |{\rm
cum}(w_{uj_1},R_{uj_2},w_{uj_3},R_{uj_4})|&=&O(|j_2|^{-1/2}|j_4|^{-1/2})\\
\label{eq:rrrr4}
|\E(w_{uj_1}R_{uj_2}w_{uj_3}R_{uj_4})|&=&O(|j_2|^{-1/2}|j_4|^{-1/2})
\end{eqnarray*}
 hold
uniformly in
$j_1,j_2,j_3,j_4=-[n/2],\cdots,-2,-1,1,2,\cdots,[n/2]$.
\end{lemma}


\noindent Proof of Lemma~\ref{lem:bound1}: Let
$K_n(\lambda)=|D_n(\lambda)|^2/(2\pi n)$ be Fej\'{e}r's kernel. We
shall first state a useful fact:
\begin{eqnarray}
\label{eq:Pj}
P_j:=\int_{-\pi}^{\pi}\left|\frac{A(\lambda)}{A(\lambda_j)}-1\right|K_n(\lambda-\lambda_j)d\lambda=O(|j|^{-1}),~\mbox{as}~n\rightarrow\infty
\end{eqnarray}
holds uniformly in $j=-[n/2],\cdots,-2,-1,1,2,\cdots,[n/2]$. The
proof of (\ref{eq:Pj}) is basically a repetition of the argument
in Robinson's (1995) Lemma 3 and is skipped. Note that
\begin{eqnarray*}
&&\hspace{-0.6cm}\E(R_{uj_1}R_{uj_2}R_{uj_3}R_{uj_4})=\cum(R_{uj_1},R_{uj_2},R_{uj_3},R_{uj_4})+\cov(R_{uj_1},R_{uj_2})\times \\
&&\hspace{-0.5cm}\cov(R_{uj_3},R_{uj_4})+\cov(R_{uj_1},R_{uj_3})\cov(R_{uj_2},R_{uj_4})+\cov(R_{uj_1},R_{uj_4})\cov(R_{uj_2},R_{uj_3})
\end{eqnarray*}
and
\begin{eqnarray*}
&&\hspace{-0.6cm}\E(w_{uj_1}R_{uj_2}w_{uj_3}R_{uj_4})=\cum(w_{uj_1},R_{uj_2},w_{uj_3},R_{uj_4})+\cov(w_{uj_1},R_{uj_2})\times\\
&&\hspace{-0.5cm}\cov(w_{uj_3},R_{uj_4})+\cov(w_{uj_1},w_{uj_3})\cov(R_{uj_2},R_{uj_4})+\cov(w_{uj_1},R_{uj_4})\cov(R_{uj_2},w_{uj_3}).
\end{eqnarray*}
 Let $T_j(\lambda)=[A(-\lambda)-A_j]|A_j|^{-1}$. We shall find a bound for each term on the RHS of the equations above.  After
some straightforward calculations, we have
\begin{eqnarray*}
|\cov(w_{uj_1},R_{uj_2})|&=&\left|\frac{1}{\sqrt{2\pi}
n}\int_{-\pi}^{\pi}
D(\lambda_{j_1}-\lambda)D(\lambda+\lambda_{j_2})T_{j_2}(\lambda) d\lambda\right|\\
&\le&\sqrt{2\pi}\left|\int_{-\pi}^{\pi}K_n(\lambda_{j_1}-\lambda)d\lambda\right|^{1/2}
P_{j_2}^{1/2}=O(j_2^{-1/2})~~\mbox{and}\\
 |\cov(R_{u j_1},{R_{u
j_2}})|&=&\left|\frac{1}{n }\int_{-\pi}^{\pi}T_{j_1}(\lambda)
T_{j_2}(-\lambda)
D_n(\lambda+\lambda_{j_1})D_n(-\lambda+\lambda_{j_2})d\lambda\right|\\
&\le& 2\pi P_{j_1}^{1/2}P_{j_2}^{1/2}=O(j_1^{-1/2} j_2^{-1/2}),
\end{eqnarray*}
where we have applied the Cauchy-Schwarz inequality. Let
$\Pi_3=[-\pi,\pi]^3$. By the Cauchy-Schwarz inequality and the
periodicity,
\begin{eqnarray*}
&&\hspace{-0.6cm}|\cum(R_{uj_1},R_{uj_2},R_{uj_3},R_{uj_4})|=\frac{|c_4(u)|}{n^2}\left|\int_{\Pi^3}D_n(\lambda_{j_1}-\lambda_1-\lambda_2-\lambda_3)D_n(\lambda_1+\lambda_{j_2})\right.\\
&&\hspace{-0.5cm}\left.D_n(\lambda_2+\lambda_{j_3})D_n(\lambda_3+\lambda_{j_4})T_{j_1}(-(\lambda_1+\lambda_2+\lambda_3)) T_{j_2}(\lambda_1)T_{j_3}(\lambda_2)T_{j_4}(\lambda_3)d\lambda_1d\lambda_2d\lambda_3\right|\\
&&\le
C\left[\int_{\Pi^3}K_n(\lambda_{j_1}-\lambda_1-\lambda_2-\lambda_3)|T_{j_1}(-(\lambda_1+\lambda_2+\lambda_3))|^2
d\lambda_1
d\lambda_2 d\lambda_3\right]^{1/2}\times\\
&&\left[
\int_{\Pi^3}K_n(\lambda_1+\lambda_{j_2})K_n(\lambda_2+\lambda_{j_3})
K_n(\lambda_3+\lambda_{j_4})|T_{j_2}(\lambda_1)|^2 |T_{j_3}(\lambda_2)|^2\times\right.\\
&&\left.|T_{j_4}(\lambda_3)|^2 d\lambda_1 d\lambda_2
d\lambda_3\right]^{1/2}\le C P_{j_1}^{1/2}
P_{j_2}^{1/2}P_{j_3}^{1/2}
P_{j_4}^{1/2}=O\left(\Pi_{i=1}^{4}j_i^{-1/2}\right).
\end{eqnarray*}
Similarly,
\begin{eqnarray*}
&&|\cum(w_{uj_1},R_{uj_2},w_{uj_3},R_{uj_4})|=\left|\frac{c_4(u)}{2\pi n^2}\int_{\Pi^3}D_n(-(\lambda_1+\lambda_2+\lambda_3))D_n(\lambda_{j_2}+\lambda_1)\right.\\
&&\hspace{0.5cm}\left.D_n(\lambda_2)D_n(\lambda_{j_4}+\lambda_3)T_{j_2}(\lambda_1)T_{j_4}(\lambda_3)d\lambda_1d\lambda_2
d\lambda_3\right|\\
&&\le C\left[\int_{\Pi^3}K_n(\lambda_{j_2}+\lambda_1)K_n(\lambda_1)K_n(\lambda_{j_4}+\lambda_3)|T_{j_2}(\lambda_1)|^2|T_{j_4}(\lambda_3)|^2d\lambda_1d\lambda_2d\lambda_3\right]^{1/2}\\
&&\hspace{0.5cm}\times\left[\int_{\Pi^3}K_n(\lambda_1+\lambda_2+\lambda_3)d\lambda_1
d\lambda_2 d\lambda_3\right]^{1/2}\le C P_{j_2}^{1/2}
P_{j_4}^{1/2}=O(j_2^{-1/2} j_4^{-1/2}).
\end{eqnarray*}
Therefore the conclusion follows. \qed

\newpage

\begin{center}
\begin{tabular}{cccccc|ccc}
\hline \hline
(a)&&&&$T_n(\hat{\pmb{\theta}})$&&&$T_n(\hat{\pmb{\gamma}})$&\\
\hline
$n$&$B_n$&$\alpha\%$&\mbox{BAR}&\mbox{TUK}&\mbox{PAR}&\mbox{BAR}&\mbox{TUK}&\mbox{PAR}\\
 \hline
64&6&$5\%$& 8.68 &8.00& 8.44&8.16 &7.58& 8.38
\\
&&$10\%$& 13.46& 12.14 &12.88 &12.62& 11.76& 12.10
\\
&10&$5\%$& 7.28 &6.94& 7.40& 6.42& 6.00& 7.00 \\
&&$10\%$&12.34& 11.50& 12.04 &11.40& 11.00& 11.14
\\
&15&$5\%$&6.86& 6.66& 6.80& 5.80& 5.44& 6.04
 \\
&&$10\%$&11.84& 11.30& 11.50& 10.40& 10.16& 10.80
\\
\hline
128&7&$5\%$&7.78& 7.16& 7.58& 7.44& 6.96&  7.7 \\
&&$10\%$&12.30& 11.68 &11.94& 11.68 &11.22 &11.28
\\
&12&$5\%$& 7.28& 7.08& 7.16& 6.52& 5.88&  6.4
\\
&&$10\%$&11.40& 11.12 &11.20& 10.84& 10.18& 10.46 \\
&20&$5\%$& 6.82 &6.70 &6.98 &5.86 &5.50  &5.8
\\
&&$10\%$&11.38& 11.04 &11.26& 10.30&  9.88 &10.24
\\
\hline\hline
\end{tabular}
\end{center}
\begin{table}
 \caption{Rejection rates in percentage under the null hypothesis: (a) when the data are generated from
model~(\ref{eq:ar1}). (b) when the data are generated from
model~(\ref{eq:ma1}). }
\label{tb:table1}
\end{table}

\begin{center}
\begin{tabular}{cccccc|ccc}
\hline \hline
(b)&&&&$T_n(\hat{\pmb{\theta}})$&&&$T_n(\hat{\pmb{\gamma}})$&\\
\hline
$n$&$B_n$&$\alpha\%$&\mbox{BAR}&\mbox{TUK}&\mbox{PAR}&\mbox{BAR}&\mbox{TUK}&\mbox{PAR}\\
 \hline
64&6&$5\%$&11.22 &10.64&  9.98& 7.36& 7.02& 7.22
\\
&&$10\%$&15.82& 14.92& 14.40& 11.70& 10.82& 11.06
\\
&10&$5\%$& 11.86& 11.66& 10.86&  7.12& 6.72& 6.94 \\
&&$10\%$&17.26& 17.02& 15.46& 11.68 &11.44& 11.50
\\
&15&$5\%$&12.66& 12.68& 11.78 & 6.68& 6.72& 6.82
 \\
&&$10\%$&18.36 &17.88& 17.28 &11.80& 11.28& 11.52
\\
\hline 128&7&$5\%$& 11.40& 10.80 &10.50& 6.64& 6.48& 6.60
 \\
&&$10\%$& 16.18 &15.48& 14.56& 10.90& 10.40 &10.34
\\
&12&$5\%$&12.50 &12.00 &11.48& 6.36& 6.38& 6.20 \\
&&$10\%$&17.30 &17.14& 16.04& 10.34& 10.44& 10.48
 \\
&20&$5\%$&13.52 &13.22& 12.48& 5.38& 5.42& 5.98
\\
&&$10\%$&19.24& 18.96 &17.58& 10.22 &10.26& 10.34
\\
\hline\hline
\end{tabular}
\end{center}

\newpage

\begin{center}
\begin{tabular}{cccccc|ccc}
\hline \hline
(a)&&&&$T_n(\hat{\pmb{\theta}})$&&&$T_n(\hat{\pmb{\gamma}})$&\\
\hline
$n$&$B_n$&$\alpha\%$&\mbox{BAR}&\mbox{TUK}&\mbox{PAR}&\mbox{BAR}&\mbox{TUK}&\mbox{PAR}\\
\hline 64&6&$5\%$& 84.00 &81.38& 84.44& 82.72& 79.28& 83.02
\\
&&$10\%$&89.56& 87.22& 89.78& 88.62& 86.22 &89.00 \\
&10&$5\%$&  79.04& 73.04 &79.34& 77.06& 72.30& 77.32\\
&&$10\%$&85.62& 82.26& 86.08& 84.66& 80.38& 85.52
\\
&15&$5\%$&72.46& 65.80& 72.58& 70.20& 63.64& 71.80
\\
&&$10\%$&80.44 &75.94& 81.56& 80.00 &74.04& 79.66
\\
\hline 128&7&$5\%$&99.14 &98.56& 99.28& 98.90& 98.38& 98.96
\\
&&$10\%$&99.58 &99.36& 99.58& 99.54 &99.20& 99.56
\\
&12&$5\%$&98.00& 96.80& 98.06& 97.94 &96.60& 97.84
\\
&&$10\%$&99.00& 98.40& 99.16& 98.90& 98.44& 99.04
\\
&20&$5\%$&95.90& 93.12& 96.12& 95.50& 92.70& 95.70
\\
&&$10\%$& 98.00& 96.58& 98.02& 97.60& 96.22 &97.80
\\
\hline\hline
\end{tabular}
\end{center}

\begin{table}
 \caption{Rejection rates in percentage under
Alternative  1: (a) when the data are generated from
model~(\ref{eq:ar1}). (b) when the data are generated from
model~(\ref{eq:ma1}). }
\label{tb:table2}
\end{table}


\begin{center}
\begin{tabular}{cccccc|ccc}
\hline \hline
(a)&&&&$T_n(\hat{\pmb{\theta}})$&&&$T_n(\hat{\pmb{\gamma}})$&\\
\hline
$n$&$B_n$&$\alpha\%$&\mbox{BAR}&\mbox{TUK}&\mbox{PAR}&\mbox{BAR}&\mbox{TUK}&\mbox{PAR}\\
\hline 64&6&$5\%$&82.06& 77.92& 83.28& 85.50& 81.66& 86.48
\\
&&$10\%$& 89.02& 86.46& 90.40& 90.90& 88.76& 91.52\\
&10&$5\%$&73.92& 68.04& 74.48& 78.90 &72.94 &79.24\\
&&$10\%$&83.20& 79.34& 84.40& 86.16 &81.64& 86.60
\\
&15&$5\%$& 65.72& 58.20& 66.54& 71.44& 63.02& 71.64
\\
&&$10\%$& 76.70& 71.10& 78.34&  80.60& 74.80& 81.14
\\
\hline 128&7&$5\%$&98.84 &97.80& 99.02 &99.42& 99.10& 99.52
\\
&&$10\%$&99.50& 99.24& 99.66& 99.76& 99.60& 99.76
\\
&12&$5\%$&97.20& 95.20& 97.24& 98.72& 97.62 &98.86
\\
&&$10\%$&98.70& 97.90& 98.96& 99.44& 98.94& 99.44
\\
&20&$5\%$& 93.50& 90.02& 93.74& 96.96& 94.42& 97.02
\\
&&$10\%$& 96.98& 94.86& 97.38& 98.38& 97.12& 98.66
\\
\hline\hline
\end{tabular}
\end{center}

\newpage

\begin{center}
\begin{tabular}{cccccc|ccc}
\hline \hline
(a)&&&&$T_n(\hat{\pmb{\theta}})$&&&$T_n(\hat{\pmb{\gamma}})$&\\
\hline
$n$&$B_n$&$\alpha\%$&\mbox{BAR}&\mbox{TUK}&\mbox{PAR}&\mbox{BAR}&\mbox{TUK}&\mbox{PAR}\\
\hline 64&6&$5\%$&61.48 &63.48& 57.98& 58.28& 59.88& 55.26
\\
&&$10\%$&73.80& 73.98& 70.12& 70.68 &71.54& 67.90 \\
&10&$5\%$& 68.04& 66.98& 65.06& 64.22& 64.08& 61.58\\
&&$10\%$&78.56& 78.48& 76.42& 75.12& 74.80& 73.58
\\
&15&$5\%$&71.40& 69.62& 68.16& 66.56& 66.04& 65.20
\\
&&$10\%$& 81.24 &80.76& 79.10& 79.00& 77.32& 75.52
\\
\hline 128&7&$5\%$&95.58 &95.36& 92.76 &93.60& 93.68& 91.48
\\
&&$10\%$&97.98& 98.08& 96.34& 97.12& 97.14& 95.38
\\
&12&$5\%$&98.52& 98.02& 96.94& 98.02& 97.32& 95.42
\\
&&$10\%$&99.52& 99.32& 98.54 &99.20& 98.92& 98.04
\\
&20&$5\%$& 99.36& 99.30& 98.90& 98.82& 98.98& 98.38
\\
&&$10\%$& 99.80& 99.82& 99.74& 99.64& 99.66 &99.42
\\
\hline\hline
\end{tabular}
\end{center}
\begin{table}
 \caption{Rejection rates in percentage under
Alternative 2: (a) when the data are generated from
model~(\ref{eq:ar1}). (b) when the data are generated from
model~(\ref{eq:ma1}).}
 \label{tb:table3}
\end{table}

\begin{center}
\begin{tabular}{cccccc|ccc}
\hline \hline
(a)&&&&$T_n(\hat{\pmb{\theta}})$&&&$T_n(\hat{\pmb{\gamma}})$&\\
\hline
$n$&$B_n$&$\alpha\%$&\mbox{BAR}&\mbox{TUK}&\mbox{PAR}&\mbox{BAR}&\mbox{TUK}&\mbox{PAR}\\
\hline 64&6&$5\%$&78.10& 76.14& 76.40& 77.82& 75.46& 75.12
\\
&&$10\%$&86.44 &84.76& 84.54 &85.98& 84.72& 83.86  \\
&10&$5\%$&76.50& 72.84& 75.60& 77.58& 73.60& 75.92\\
&&$10\%$&86.80& 84.08 &85.20&  86.78& 84.76 &85.58
\\
&15&$5\%$&75.62& 71.54& 73.44&  76.80 &72.06& 74.22
\\
&&$10\%$&86.98& 83.38& 84.22& 87.18 &84.00& 85.28
\\
\hline 128&7&$5\%$&97.76 &97.52&  100 &99.20& 98.94& 98.64
\\
&&$10\%$&99.46& 99.40& 99.16& 99.68 &99.60& 99.48
\\
&12&$5\%$&98.44& 97.34&  100 &99.50& 99.12& 99.04
\\
&&$10\%$&99.54& 99.24& 99.34& 99.80& 99.72& 99.68
\\
&20&$5\%$& 98.48& 98.06&  100 & 99.68& 99.58& 99.44
\\
&&$10\%$& 99.62 &99.48& 99.40& 99.84 &99.84& 99.78
\\
\hline\hline
\end{tabular}
\end{center}

\newpage

\begin{center}
\begin{tabular}{cccccc|ccc}
\hline \hline
(a)&&&&$T_n(\hat{\pmb{\theta}})$&&&$T_n(\hat{\pmb{\gamma}})$&\\
\hline
$n$&$B_n$&$\alpha\%$&\mbox{BAR}&\mbox{TUK}&\mbox{PAR}&\mbox{BAR}&\mbox{TUK}&\mbox{PAR}\\
\hline 64&6&$5\%$&17.98 &16.64&  6.74& 13.34& 11.62&  4.28
\\
&&$10\%$& 31.48& 27.30& 13.26&  25.8& 22.62&  9.48 \\
&10&$5\%$&37.06& 42.36& 27.24& 33.22& 40.26& 22.10\\
&&$10\%$&51.16& 57.38& 41.38& 47.6& 53.58& 37.88
\\
&15&$5\%$&42.30 &46.98& 43.32& 37.90& 44.36& 41.28
\\
&&$10\%$&55.54& 60.26& 58.16& 53.1 &57.28& 54.12
\\
\hline 128&7&$5\%$&71.62& 70.66& 22.38& 65.22& 65.08& 13.44
\\
&&$10\%$&83.18& 83.00& 37.24& 80.44& 80.70& 28.02
\\
&12&$5\%$&85.76& 89.52& 82.76& 84.68& 89.16 &80.72
\\
&&$10\%$&92.78& 94.44 &91.50& 91.52& 94.10 &90.40
\\
&20&$5\%$& 86.80& 88.84& 89.66& 85.00& 87.22& 89.04
\\
&&$10\%$& 92.60 &93.50& 94.36& 91.92& 93.12& 93.80
\\
\hline\hline
\end{tabular}
\end{center}
\begin{table}
 \caption{Rejection rates in percentage under
Alternative 3: (a) when the data are generated from
model~(\ref{eq:ar1}). (b) when the data are generated from
model~(\ref{eq:ma1}).}
 \label{tb:table4}
\end{table}

\begin{center}
\begin{tabular}{cccccc|ccc}
\hline \hline
(a)&&&&$T_n(\hat{\pmb{\theta}})$&&&$T_n(\hat{\pmb{\gamma}})$&\\
\hline
$n$&$B_n$&$\alpha\%$&\mbox{BAR}&\mbox{TUK}&\mbox{PAR}&\mbox{BAR}&\mbox{TUK}&\mbox{PAR}\\
\hline 64&6&$5\%$& 17.38 &14.86&  6.46 &13.86 &11.08&  3.52
\\
&&$10\%$&32.18& 28.62& 14.04& 28.1& 23.16 & 8.10 \\
&10&$5\%$& 33.10 &38.64& 23.14&  33.44 &40.12& 22.60\\
&&$10\%$&50.18 &55.18& 40.86& 49.6& 55.48& 37.98
\\
&15&$5\%$&38.20& 42.52& 38.90& 39.38 &44.38& 40.46
\\
&&$10\%$& 52.96 &56.98& 55.48& 54.6& 60.26& 56.56
\\
\hline 128&7&$5\%$&63.60& 64.06& 22.40& 70.72& 69.52& 13.00
\\
&&$10\%$& 80.98 &80.50& 39.56&  84.36& 83.54 &28.66
\\
&12&$5\%$&81.02& 85.44& 77.04& 86.74& 90.16& 84.06
\\
&&$10\%$& 90.24& 92.64& 88.92& 93.02& 94.70& 91.56
\\
&20&$5\%$& 81.54& 83.86& 85.44&  87.58& 89.08& 90.18
\\
&&$10\%$& 89.80 &91.46& 92.44 &92.90& 93.88& 94.52
\\
\hline\hline
\end{tabular}
\end{center}

\baselineskip=12pt
\bigskip
\end{document}